\documentclass[11pt,psamsfonts]{amsart}
\usepackage{amsmath}
\usepackage{amsthm}
\usepackage{amssymb}
\usepackage{amscd}
\usepackage{amsfonts}
\usepackage{amsbsy}
\usepackage[latin1]{inputenc}
\usepackage{epsfig,afterpage}
\usepackage[dvips]{psfrag}
\usepackage{psfrag}
\usepackage[all]{xy}
\usepackage{color} 
\usepackage{multirow}
\usepackage{pdflscape}

\usepackage{overpic}

\newtheorem {theorem} {Theorem}

\begin{document}

\title[Periodic Solutions of the $\phi^{6}$-Van der Pol Oscillator] {Periodic Solutions of the extended Duffing-Van der Pol Oscillator}

\author[R.D. Euzébio and J. Llibre]
{Rodrigo D. Euzébio$^{1,2}$ and Jaume Llibre$^1$}

\address{$^1$ Departament de Matem\`{a}tiques,
Universitat Aut\`{o}noma de Barcelona, 08193 Bellaterra, Barcelona,
Catalonia, Spain} \email{jllibre@mat.uab.cat}

\address{$^2$ Departament de Matem\'atica,
IBILCE, UNESP, Rua Cristovao Colombo 2265, Jardim Nazareth, CEP
15.054--00, Sao Jos\'e de Rio Preto, SP, Brazil}
\email{rodrigo.euzebio@sjrp.unesp.br}



\date{}
\dedicatory{}

\begin{abstract}
In this paper some aspects on the periodic solutions of the extended Duffing-Van der Pol oscillator are discussed. Doing different rescaling of the variables and parameters of the system associated to the extended Duffing-Van der Pol oscillator we show that it can bifurcate one or three periodic solutions from a 2--dimensional manifold filled by periodic solutions of the referred system. For each rescaling we exhibit concrete values for which these bounds are reached. Beyond that we characterize the stability of some periodic solutions. Our approach is analytical and the results are obtained using averaging theory and some algebraic techniques.
\end{abstract}

\maketitle

\section{Introduction}\label{secao introducao}

\subsection{Setting the problem}\label{secao colocacao do problema}

A large number of non--autonomous chaotic phenomena in physics, engineering, mechanics, biology among others are described by second-order differential systems of the form
\begin{equation}\label{general_system}
\ddot{x}=g(x,\dot{x},t)+\gamma(t),
\end{equation}
where $g(x,\dot{x},t)$ is a continuous function and $\gamma(t)$ is some external force. For instance, in biology system \eqref{general_system} models the FHN neuron oscillator and in the engineering system \eqref{general_system} it is a model to the horizontal platform system.

\smallskip

The specific topic addressed in this paper concerns with other particular case of equation \eqref{general_system}, namely, an extension of the forced Van der Pol equation with external excitation. Van der Pol's system also plays an important role in many applications in areas like engineering, biology, physics and seismology (see \cite{GUCKII} and references therein).

\smallskip

The system associated with the forced Van der Pol equation with external excitation is characterized by system \eqref{general_system} with functions $g$ and $\gamma$ in the form
\begin{equation}\label{special_functions}
g(x,\dot{x},t)=\rho_{0}(1-x^2)\dot{x}-\dfrac{d}{dx}V(x),\quad \gamma(t)=\delta_{0}\cos(\omega t),
\end{equation}
where $\rho_{0}$ is the damping parameter, $V(x)$ is the potential function and $\delta_{0}$ and $\omega$ are the amplitude and angular frequency of the driving force $\gamma(t)$, respectively. We assume that $\rho_{0}$ is non--negative and $\delta_{0}$ and $\omega$ are positive. The potential $V(x)$ can be approximated by a finite Taylor expansion in series
$$
V_{2}(x)=\dfrac{1}{2}\omega_{0}^2x^2,\;\; V_4(x)=\dfrac{1}{2}\omega_{0}^2x^2+\dfrac{1}{4}\alpha_{0} x^4,\;\; V_6(x)=\dfrac{1}{2}\omega_{0}^2x^2+\dfrac{1}{4}\alpha_{0} x^4+\dfrac{1}{6}\lambda_{0} x^6,
$$
where $\omega_0$ and $\lambda_{0}$ are non--zero and $\alpha_{0}$ is a real number.

\smallskip

Almost all papers on forced excited Van der Pol systems deal with the potential $V_2$. However some papers concerning the potential $V_{4}$ have shown a lot of interesting behaviours (see \cite{CJN}, \cite{DV7}, \cite{DV5}, \cite{DV4} and \cite{DV6}). This case is usually refereed as Duffing--Van der Pol oscillator or $\phi^{4}-$Van der Pol oscillator. Nevertheless more recently some papers were published taking into account the potential $V_{6}$ meanly addressing with the problem of chaos control (see \cite{PRINCIPAL20}, \cite{PRINCIPAL} and \cite{PRINCIPAL21}). The dynamics considering the potential $V_{6}$ is more complex and rich than the corresponding cases considering the potentials $V_{2}$ and $V_{4}$ (see \cite{COMPLEX9} and \cite{COMPLEX8}). This is the case that we will regard. This case is quoted in the literature as extended Duffing-Van der Pol oscillator or $\phi^{6}-$Van der Pol oscillator.

\smallskip

In this paper we will give an analytical treatment to system \eqref{general_system} in order to study its periodic solutions considering functions \eqref{special_functions} and the potential $V_{6}$. Indeed, calling $y=\dot{x}$ we obtain
\begin{equation}
\begin{array}{l}
\dot{x}=y,\vspace{0.2cm}\\
\dot{y}=-\omega_{0}^{2}x+\rho_{0} y-\alpha_{0} x^3-\rho_{0} x^2 y-\lambda_{0} x^5+\delta_{0}\cos(\omega t).
\end{array}\label{vdp_original}
\end{equation}

System \eqref{vdp_original} becomes simpler if we perform a rescaling $s=\omega t$ in the time $t$ and another one $y=\omega_{0} Y$ in the spatial variable $y$. In fact, calling again the new time $s$ by $t$ and the variable $Y$ by $y$, after the rescaling we have
\begin{equation}\label{vdp}
\begin{array}{l}
\dot{x}=y,\vspace{0.2cm}\\
\dot{y}=-x+\rho y-\alpha x^3-\rho x^2 y-\lambda x^5+\delta\cos\,t,
\end{array}
\end{equation}
where the new parameters $\rho$, $\alpha$, $\lambda$ and $\delta$ are the respective old ones divided by $\omega_{0}^{2}$. The study of the periodic solutions of the non--linear non--autonomous $2\pi-$periodic differential system \eqref{vdp} will be the objective of this paper.

\smallskip

When $\alpha$, $\lambda$ and $\delta$ are zero equation $\eqref{vdp}$ is referred as unforced Van der Pol equation and has a unique stable periodic solution for $\rho$ positive. Furthermore, if $\rho$ is large this periodic solution remains and it describes a periodic oscillatory behaviour called relaxation oscillation. On the other hand, by considering $\delta$ non--zero, system $\eqref{vdp}$ can present nonlinear attractors from simple periodic solutions up to multiperiodicity and chaos. In \cite{DIDATIC} we can find a didactical discussion about these objects and some results using classical techniques are obtained for a special case of system \eqref{vdp} using particular values of $\rho$, $\delta$, $\omega$ and initial conditions.

\smallskip

There exists an exhaustive list of papers in the literature studying the properties of system \eqref{vdp} when $\delta$ is zero. For $\delta$ positive many open questions remain mainly due to the difficulty of integrate the system. In particular, in \cite{GENERAL} the authors study some aspects of robust practical synchronization for system \eqref{general_system} and apply the results to a particular case of system \eqref{vdp}. In the same direction, in \cite{LEUNG} the author investigate synchronization processes between chaotic attractors considering $\alpha=\lambda=0$ in system \eqref{vdp}. In \cite{GUCKII} and \cite{GUCKI} it is studied many aspects of system \eqref{vdp} for the case where $\alpha$ and $\lambda$ are zero and $\rho$ is large. Furthermore in \cite{EH} the authors provide sufficient conditions to the existence of one periodic solution under some analytical hypotheses for a general forced Van der Pol system considering $\rho$ equal to zero.

\smallskip

Periodic solutions for system \eqref{vdp} were found in \cite{PRINCIPAL}, where the authors investigate chaos control. Besides, in \cite{CJ} and \cite{CL} were obtained results on the existence of periodic solutions for an autonomous special case of the extended Duffing-Van der Pol oscillator. In \cite{LXZ}, \cite{YXS} and \cite{YXS2} we also can find some results on systems similar to system \eqref{vdp}. 

\smallskip

In this paper we are concerned with periodic solutions of system \eqref{vdp}. We will present sufficient conditions in order that this system possess one or three periodic solution and we will provide conditions on the parameters $\rho$, $\alpha$, $\lambda$ and $\delta$ for which these bounds are realizable. In Addition
 we will prove that it is not possible to obtain different bounds of periodic solutions using the methodology presented in this paper.

\smallskip

The phase space of our non--autonomous differential system is $(x,y,t)$ and its Poincaré map is defined in the $(x,y)$--space. When possible the stability of the periodic solutions will be studied. As usual a periodic solution is {\it stable} if the eigenvalues of the fixed point associated to its Poincaré map have negative real part, otherwise there periodic solution is \textit{unstable}. Inside the unstable periodic solutions the are two types. The unstable saddle periodic solutions having an eigenvalue with negative real part and the other with positive real part, and the repellor with both eigenvalues having negative real part.



\smallskip

We should note that the method used here for studying periodic solutions can be applied to any periodic non--autonomous differential system as done in \cite{EL} and \cite{LV}. In these papers the authors applied the method used in this paper in order to guarantee the existence of periodic solutions in a periodic FitzHugh--Nagumo system and in the Vallis system, respectively.

\smallskip

The paper is organized as follows. In subsection \ref{secao enunciados} the main results are stated and compared with other results. Also, some important points on the results are clarified. In subsection \ref{proofs}  we prove the results. In subsection \ref{discussion} some aspects on the results are pointed out. Lastly, section \ref{averaging} is devoted to give a brief summary of the results that we use from averaging theory.

\subsection{Statement of the main results}\label{secao enunciados}

In this section we present our results. We will behave under an analytical approach and for this reason the results are valid for a large range of the parameters of system \eqref{vdp}, different from the major part of the results dealing with numerical techniques. In addition, it is important to note that we are concerned with harmonic solutions in the sense that they do not bifurcate from periodic solutions with multiple period.

%

\smallskip

We have the following results.

\smallskip


\begin{theorem}\label{theoremc8}
Consider $\varepsilon>0$ sufficiently small, $(\rho,\delta,\alpha,\lambda)=(\varepsilon r,\varepsilon d,\varepsilon^{n_{2}}a,$  $\varepsilon^{n_{3}}\ell)$ with $n_{2},n_{3}>1$ and $81d^{2}>48r^{2}>0$. Then system \eqref{vdp} has a $2\pi$--periodic solution $(x(t,\varepsilon),y(t,\varepsilon))$ such that
$$
(x(0,\varepsilon),y(0,\varepsilon))\to\left( 0,\dfrac{2(36r)^{1/3}}{3\Gamma}+\dfrac{1}{3}\left(\dfrac{6}{r}\right)^{1/3}\Gamma   \right),
$$
when $\varepsilon\to 0$, where $\Gamma=\left(9d+\sqrt{81d^{2}-48r^{2}}\right)^{1/3}$. Moreover, for $r$ sufficiently small this periodic solution is stable.
\end{theorem}

\smallskip


\begin{theorem}\label{theoremc10}
Consider $\varepsilon>0$ sufficiently small, $(\alpha,\delta,\rho,\lambda)=(\varepsilon a,\varepsilon d,\varepsilon^{n_{1}} r,$ $\varepsilon^{n_{3}} \ell)$ with $n_{1},n_{3}>1$ and $\alpha\neq0$. Then system \eqref{vdp} has a $2\pi$--periodic solution $(x(t,\varepsilon),y(t,\varepsilon))$ such that
$$
(x(0,\varepsilon),y(0,\varepsilon))\to\left(\left(   \dfrac{4d}{3a}  \right)^{1/3} ,0   \right),
$$
when $\varepsilon\to 0$.
\end{theorem}

\smallskip


\begin{theorem}\label{theoremc11}
Consider $\varepsilon>0$ sufficiently small and $(\lambda,\delta,\rho,\alpha)=(\varepsilon \ell,\varepsilon d,\varepsilon^{n_{1}} r,$ $\varepsilon^{n_{2}} a)$ with $n_{1},n_{2}>1$. Then system \eqref{vdp} has a $2\pi$--periodic solution $(x(t,\varepsilon),$ $y(t,\varepsilon))$ such that
$$
(x(0,\varepsilon),y(0,\varepsilon))\to\left(\left(   \dfrac{8d}{5\ell}  \right)^{1/5},0  \right),
$$
when $\varepsilon\to 0$.
\end{theorem}

\smallskip

We remark that the stability of the periodic solutions of Theorems 2 and 3 cannot be decided with the real part of the eigenvalues of their Poincaré map because these real part are zero.

\smallskip


\begin{theorem}\label{theoremc17}
Consider $\varepsilon>0$ sufficiently small and $(x,y,\rho,\delta,\alpha,\lambda)=(\varepsilon X,$ $\varepsilon Y,$ $\varepsilon r,\varepsilon^{2} d,\varepsilon^{n_{2}}a,\varepsilon^{n_{3}}\ell)$ with $\rho\neq0$. Then system \eqref{vdp} has a $2\pi$--periodic solution $(x(t,\varepsilon),y(t,\varepsilon))$ such that
$$
(x(0,\varepsilon),y(0,\varepsilon))=\left( \mathit{O}(\varepsilon),-\dfrac{\delta}{\varepsilon\rho}+\mathit{O}(\varepsilon)   \right),
$$
when $\varepsilon\to 0$. Moreover this periodic solution is unstable.
\end{theorem}

\smallskip

Note that in Theorem 4 we have a periodic solution that comes from infinity. As far as we know this kind of behaviour also has not been observed in the papers concerned with system \eqref{vdp}. Nevertheless this behaviour is common when we perform a rescaling in the parameters and variables as we made previously and can be observed also in papers \cite{EL} and \cite{LV}.

\smallskip


\begin{theorem}\label{theoremc5}
Consider $\varepsilon>0$ sufficiently small, $(\alpha,\lambda,\delta,\rho)=(\varepsilon a,\varepsilon \ell,\varepsilon d,$ $\varepsilon^{n_{1}} r)$ with $n_{1}>1$, $\alpha\lambda<0$ and
$$
D=\dfrac{53747712 a^5 d^2}{78125 \ell^7}+\dfrac{20480 d^4}{\ell^{4}}.
$$
Then system \eqref{vdp} has a $2\pi$--periodic solution if $D=0$ and three $2\pi$--periodic solution if $D<0$. Moreover, there are values of $\alpha$, $\lambda$ and $\delta$ that realize these number of periodic solutions.
\end{theorem}

\smallskip


\begin{theorem}\label{theoremc3}
Consider $\varepsilon>0$ sufficiently small, $(\rho,\alpha,\delta,\lambda)=(\varepsilon r,\varepsilon a,\varepsilon d,$ $\varepsilon^{n_{3}}\ell)$ with $n_{3}>1$ and $-324 a^4 r^2-9 a^2 d^2 r^2+36 a^2 r^4-d^2 r^4\neq0$. Consider also the values
$$
\begin{array}{l}
\Delta_{1}=324 a^4 + d^2 r^2 + 9 a^2(d^2 - 4 r^2),\vspace{0.2cm}\\
\Delta_{2}=2187a^4d^4+27d^4r^4-16d^2r^6+18a^2(27d^4r^2-72d^2r^4+32r^6).
\end{array}
$$
The system
$$
\begin{array}{r}
rx_{0}(-4+x_{0}^{2}+y_{0}^{2})-3ay_{0}(x_{0}^{2}+y_{0}^{2})=0,\vspace{0.2cm}\\
4d-ry_{0}(-4+x_{0}^{2}+y_{0}^{2})-3ax_{0}(x_{0}^{2}+y_{0}^{2})=0,
\end{array}
$$
has one solution $(x_{0}^{0},y_{0}^{0})$ if $\Delta_{1}\Delta_{2}=0$ or $\Delta_{2}>0$ and three solutions $(x_{0}^{i},y_{0}^{i})$ if $\Delta_{2}<0$ for $i=1,2,3$, all of them satisfying
$$
27\alpha^{2}\left((x_{0}^{i})^{2}+(y_{0}^{i})^{2}\right)^{2}+\rho^{2}\left(  -4+(x_{0}^{i})^{2}+(y_{0}^{i})^{2}   \right)\left(  -4+3(x_{0}^{i})^{2}+3(y_{0}^{i})^{2}   \right)\neq 0.
$$
Then if $\Delta_{1}\Delta_{2}=0$ or $\Delta_{2}>0$ system \eqref{vdp} has a $2\pi$--periodic solution $(x^{0}(t,\varepsilon),$ $y^{0}(t,\varepsilon))$ such that $(x^{0}(0,\varepsilon),y^{0}(0,\varepsilon))\to(x_{0}^{0},y_{0}^{0})$ when $\varepsilon\to 0$. Additionally if $\Delta_{2}<0$ system \eqref{vdp} has three $2\pi$--periodic solution $(x^{i}(t,\varepsilon),y^{i}(t,\varepsilon))$ such that $(x^{i}(0,\varepsilon),y^{i}(0,\varepsilon))\to(x_{0}^{i},y_{0}^{i})$ when $\varepsilon\to 0$ for $i=1,2,3$. Moreover there are values of $\rho$, $\alpha$ and $\delta$ for which these number of periodic solutions is reached.
\end{theorem}

%

\smallskip


\begin{theorem}\label{theoremc1}
Consider $\varepsilon>0$ sufficiently small, $(\alpha,\rho,\lambda,\delta)=(\varepsilon a,\varepsilon r,\varepsilon\ell,\varepsilon d)$, $\alpha\rho\neq 0$, $(r^{2}-3a^{2})(r^{2}-9a^{2})\neq0$ and
$$
\begin{array}{l}
C=4 (3 a + 10 \ell) (540 a^3 \ell - 9 a^2 (d^2 - 600 \ell^2) - 
120 a \ell (d^2 - 150 \ell^2)+  \vspace{0.2cm}\\
     \quad 50 \ell^2 (-7 d^2 + 400 \ell^2)) r^6 + 
  6 (a + 5 \ell) (6 a - d + 20 \ell) (6 a + d + 20 \ell) r^8)\neq 0,\vspace{0.2cm}\\
D=2066242608a^{14}d^{6}-3125d^{8}r^{12}-531441a^{12}(3125d^{8}+96d^{6}r^{2}+\vspace{0.2cm}\\
\quad 1536d^{4}r^{4}-1024d^{2}r^{6})+354294a^{10}(3125d^{8}r^{2}+616d^{6}r^{4}-1600\vspace{0.2cm}\\
\quad d^{4}r^{6})+18a^{2}r^{10}(9375d^{8}+5000d^{6}r^{2}-88000d^{4}r^{4}+102400d^{2}r^{6}-\vspace{0.2cm}\\
\quad 32768r^{8})+2916a^{6}r^{6}(15625d^{8}+11700d^{6}r^{2}-45632d^{4}r^{4}+23552\vspace{0.2cm}\\
\quad d^{2}r^{6}-2048r^{8})-6561a^{8}r^{4}(46875d^{8}+24832d^{6}r^{2}-66816d^{4}r^{4}+\vspace{0.2cm}\\
\quad 12288d^{2}r^{6}+4096r^{8})+81a^{4}r^{8}(-46875d^{8}-36000d^{6}r^{2}+275200d^{4}\vspace{0.2cm}\\
\quad r^{4}-246784d^{2}r^{6}+61440r^{8}).
\end{array}
$$
The system
$$
\begin{array}{r}
2rx_{0}(-4+x_{0}^{2}+y_{0}^{2})-y_{0}(x_{0}^{2}+y_{0}^{2})(6a+5\ell(x_{0}^{2}+y_{0}^{2})=0,\vspace{0.2cm}\\
8d+8ry_{0}-(x_{0}^{2}+y_{0}^{2})(6ax_{0}+2ry_{0}+5\ell x_{0}(x_{0}^{2}+y_{0}^{2}))=0,
\end{array}
$$
has one solution $(x_{0}^{0},y_{0}^{0})$ if $D<0$ and three solutions $(x_{0}^{i},y_{0}^{i})$ if $D>0$ for $i=1,2,3$, all of them satisfying
$$
\begin{array}{l}
4\rho^{2}\left(-4+(x_{0}^{i})^{2}+(y_{0}^{i})^{2}\right)(-4+3(x_{0}^{i})^{2}+3(y_{0}^{i})^{2})+((x_{0}^{i})^{2}+(y_{0}^{i})^{2})^{2}(6\alpha+\vspace{0.2cm}\\
\rho\lambda((x_{0}^{i})^{2}+(y_{0}^{i})^{2}))(18\alpha+25\lambda(x_{0}^{i})^{2}+3(y_{0}^{i})^{2}))\neq0.
\end{array}
$$
Then, if $D<0$ system \eqref{vdp} has a $2\pi$--periodic solution $(x^{0}(t,\varepsilon),y^{0}(t,\varepsilon))$ such that $(x^{0}(0,\varepsilon),y^{0}(0,\varepsilon))\to(x_{0}^{0},y_{0}^{0})$ when $\varepsilon\to 0$. Furthermore, if $D>0$ system \eqref{vdp} has three $2\pi$--periodic solution $(x^{i}(t,\varepsilon),y^{i}(t,\varepsilon))$ such that $(x^{i}(0,\varepsilon),$ $y^{i}(0,\varepsilon))\to(x_{0}^{i},y_{0}^{i})$ when $\varepsilon\to 0$ for $i=1,2,3$. Moreover, there are values of $\alpha$, $\rho$, $\lambda$ and $\delta$ which realizes the number of these periodic solutions.
\end{theorem}



\smallskip


\begin{theorem}\label{theoremc4}
Consider $\varepsilon>0$ sufficiently small, $(\rho,\lambda,\delta,\alpha)=(\varepsilon r,\varepsilon\ell,\varepsilon d,$ $\varepsilon^{n_{2}} a)$ with $n_{2}>1$. Consider the numbers
$$
\begin{array}{l}
C=4000 d^3 \ell^3 r^4 + 7200000 d \ell^5 r^4 + 60 d^3 \ell r^6 + 
 72000 d \ell^3 r^6 + (50 d^4 \ell^2 r^4 -\vspace{0.2cm}\\
\quad 220000 d^2 \ell^4 r^4 - 
   8000000 \ell^6 r^4 - 2800 d^2 \ell^2 r^6 + 160000 \ell^4 r^6 - 6 d^2 r^8 + \vspace{0.2cm}\\
\quad  2400 \ell^2 r^8),\vspace{0.2cm}\\
N_{2}=\dfrac{-4500 \ell^2 r^4 + 3 r^6}{3125 d^2 \ell^4},\vspace{0.2cm}\\
N_{3}=5859375 d^4 \ell^6 + 18750000 d^2 \ell^6 r^2 + 87500 \ell^4 (-7 d^2 + 1200 \ell^2) r^4 + \vspace{0.2cm}\\
\quad 25 \ell^2 (-7 d^2 + 155600 \ell^2) r^6 + 15600 \ell^2 r^8 + 9 r^{10},\vspace{0.2cm}\\
N_{4}= 1171875 d^8 \ell^6 (1250000 \ell^4 + 8500 \ell^2 r^2 + r^4) - 2500 d^6 \ell^4 (312500000000\vspace{0.2cm}\\
\quad \ell^8 + 5343750000 \ell^6 r^2 +74625000 \ell^4 r^4 + 358375 \ell^2 r^6 + 69 r^8) -3000 \vspace{0.2cm}\\
 \quad d^4 \ell^2 r^2(218750000000 \ell^{10} - 21625000000 \ell^8 r^2 - 686625000 \ell^6 r^4 - \vspace{0.2cm}\\
\quad 3715000 \ell^4 r^6 - 650 \ell^2 r^8 - 3 r^{10}) -1600 \ell^2 r^6 (400000000000 \ell^{10}-\vspace{0.2cm}\\
\quad 144900000000 \ell^8 r^2-2817000000 \ell^6 r^4 - 8880000 \ell^4 r^6 + 38700 \ell^2 r^8 + \vspace{0.2cm}\\
\quad 27 r^{10}) +12 d^2 r^4 (-50000000000000 \ell^{12} - 32250000000000 \ell^{10} r^2-\vspace{0.2cm}\\
\quad 662400000000 \ell^8 r^4 - 2803000000 \ell^6 r^6 + 5410000 \ell^4 r^8 + 19200 \ell^2 r^{10} +9 r^{12})\vspace{0.2cm}\\
N_{5}=48828125 d^8 \ell^6 - 4687500 d^6 \ell^4 r^4 + 6400 \ell^2 r^{10} (1600 \ell^2 + r^2) -16 d^2 r^8 \vspace{0.2cm}\\
\quad (2000000 \ell^4 + 2900 \ell^2 r^2 + r^4) +  d^4 r^6 (27500000 \ell^4 + 60000 \ell^2 r^2 +27 r^4)\vspace{0.2cm}\\
M_{5}=-25d^{4}\ell{2}+d^{2}(110000\ell{4}+1400\ell{2}r^{2}+3r^{4})+400(10000\ell{6}-200\ell^{4}r^{2}-\vspace{0.2cm}\\
\quad 3\ell^{2}r^{4}),
\end{array}
$$
where $CM_{5}\neq 0$. The system
$$
\begin{array}{r}
2rx_{0}(-4+x_{0}^{2}+y_{0}^{2})-5\ell y_{0}(x_{0}^{2}+y_{0}^{2})^{2}=0,\vspace{0.2cm}\\
8d-2ry_{0}(-4+x_{0}^{2}+y_{0}^{2})-5\ell x_{0}(x_{0}^{2}+y_{0}^{2})^{2}=0,
\end{array}
$$
has one solution $(x_{0}^{0},y_{0}^{0})$ if $N_{2}\leq0$, $N_{3}\geq0$ or $N_{4}\geq0$ and $N_{5}>0$ and three solutions $(x_{0}^{i},y_{0}^{i})$ if $N_{5}<0$ for $i=1,2,3$, all of them satisfying
$$
125\lambda^{2}\left((x_{0}^{i})^{2}+(y_{0}^{i})^{4}\right)+4\rho^{2}(-4+(x_{0}^{i})^{2}+(y_{0}^{i})^{2})(-4+3(x_{0}^{i})^{2}+3(y_{0}^{i})^{2})\neq0.
$$
Hence, if $N_{2}\leq0$, $N_{3}\geq0$ or $N_{4}\geq0$ and $N_{5}>0$, then system \eqref{vdp} has a $2\pi$--periodic solution $(x^{0}(t,\varepsilon),y^{0}(t,\varepsilon))$ such that $(x^{0}(0,\varepsilon),y^{0}(0,\varepsilon))\to(x_{0}^{0},y_{0}^{0})$ when $\varepsilon\to 0$. If $N_{5}<0$ then system \eqref{vdp} has three $2\pi$--periodic solution $(x^{i}(t,\varepsilon),y^{i}(t,\varepsilon))$ such that $(x^{i}(0,\varepsilon),y^{i}(0,\varepsilon))\to(x_{0}^{i},y_{0}^{i})$ when $\varepsilon\to 0$ for $i=1,2,3$. Moreover, there are values of $\rho$, $\lambda$ and $\delta$ for which system \eqref{vdp} has one or three periodic solutions.
\end{theorem}

\smallskip

We note that in Theorems 5, 6, 7 and 8 we do not say anything about the kind of stability of the periodic solutions because we do not have the explicit expressions of the real part of the eigenvalues of their Poincaré maps.

\smallskip

We remark that the periodic solutions provided in Theorems 1 until 8 exit when we take small values for the parameters of system \eqref{vdp} obeying some relations among them. For instance Theorem 1 states the existence of one periodic solution for system \eqref{vdp} if each parameter of this system is small and $\alpha$ and $\lambda$ are much smaller than $\rho$ and $\delta$. This assertion becomes more clear if we observe the replacement done in each theorem and taking into account that $\varepsilon$ is sufficiently small. As far as we know periodic solutions of system \eqref{vdp} whose parameters have this characteristic have not been observed in the literature. Besides, it seems that the simultaneously bifurcation of three harmonic periodic solutions in extended Duffing-Van der Pol system's is also new.

%

\smallskip

Furthermore we note that the periodic solutions presented in the present paper are different from those ones stated in \cite{EH}. Indeed, in the referred paper the authors ask for $V'_{6}(0)<0$, and in our case we have $V'_{6}(0)=0$.

\section{Proof and discussion of the results}\label{proof-discussion}

\subsection{Proof of the results}\label{proofs}

In order to apply averaging theory described in section \ref{averaging} in systems \eqref{vdp} we start performing a rescaling of the variables $x$ and $y$ and of the parameters $\rho$, $\alpha$, $\lambda$ and $\delta$ as follows.
\begin{equation}\label{rescaling}
\begin{array}{c}
x=\varepsilon^{m_{1}}X,\qquad y=\varepsilon^{m_{2}}Y,\vspace{0.2cm}\\
\rho=\varepsilon^{n_{1}}r,\qquad \alpha=\varepsilon^{n_{2}}a,\qquad \lambda=\varepsilon^{n_{3}}\ell,\qquad \delta=\varepsilon^{n_{4}}d,
\end{array}
\end{equation}
where $\varepsilon$ is positive and sufficiently small and $m_{i}$ and $n_{j}$ are non--negative integers, for $i=1,2$ and $j=1,2,3,4$. We recall that  since $\delta>0$, $\rho\geq 0$ and $\lambda\neq 0$ we have $d>0$, $r\geq 0$ and $\ell\neq 0$.

\smallskip

In the new variables $(X,Y)$ system \eqref{vdp} writes
\begin{equation}\label{vdp_rescaling}
\begin{array}{rcl}
\dfrac{dX}{dt}&=&\varepsilon^{-m_{1}+m_{2}}Y,\vspace{0.2cm}\\
\dfrac{dY}{dt}&=&-\varepsilon^{m_{1}-m_{2}}X+\varepsilon^{n_{1}}rY-\varepsilon^{2m_{1}+n_{1}}rX^{2}Y-\varepsilon^{3m_{1}-m_{2}+n_{2}}aX^{3}\vspace{0.2cm}\\
&&-\varepsilon^{5m_{1}-m_{2}+n_{3}}\ell X^{5}+\varepsilon^{-m_{2}+n_{4}}d\cos\,t.
\end{array}
\end{equation}
In order to have non--negative powers of $\varepsilon$ we must impose the conditions
\begin{equation}
m_{1}=m_{2}=m,\quad\mbox{and}\quad n_{4}\geq m,\label{conditions}
\end{equation}
where $m$ is a non--negative integer.

Hence, with conditions \eqref{conditions} system \eqref{vdp_rescaling} becomes
\begin{equation}\label{vdp_new}
\begin{array}{rcl}
\dfrac{dX}{dt}&=&Y,\vspace{0.2cm}\\
\dfrac{dY}{dt}&=&-X+\varepsilon^{n_{1}}rY-\varepsilon^{2m+n_{1}}rX^{2}Y-\varepsilon^{2m+n_{2}}aX^{3}-\varepsilon^{4m+n_{3}}\ell X^{5}\vspace{0.2cm}\\
&&+\varepsilon^{-m+n_{4}}d\cos\,t.
\end{array}
\end{equation}

\smallskip

In this paper we will find periodic solutions of system \eqref{vdp_new} depending on the parameters $r$, $a$, $\ell$ and $d$ and the powers $m$, $n_{i}$ of $\varepsilon$, for $i=1,2,3,4$. Then we will go back through the rescaling \eqref{rescaling} to ensure the existence of periodic solutions in system \eqref{vdp}. 

\smallskip

From now on we assume that the values $n_{1}$, $m^2+n_{2}^2$, $m^2+n_{3}^2$ and $n_{4}-m$ are positive and observe that considering these conditions each power of $\varepsilon$ in system \eqref{vdp_new} becomes positive. The reason for a such assumption will be explained later on in section \ref{discussion}. Now we shall apply the averaging theory described in section \ref{averaging} . Thus following the notation of the mentioned section and denoting again the variables $(X,Y)$ by $(x,y)$ we have ${\bf x}=(x,y)^{T}$ and system \eqref{unperturberd} corresponding to system \eqref{vdp_new} can be written as
\begin{equation}\label{F0_vdp}
{\bf \dot{x}}=F_{0}(t,{\bf x})=(y,-x)^{T}.
\end{equation}
So the solution ${\bf x}(t,{\bf z})=(x(t,{\bf z}),y(t,{\bf z}))$ of system \eqref{F0_vdp} such that ${\bf x}(0,{\bf z})$ $={\bf z}=(x_{0},y_{0})$ is
$$
\begin{array}{rcl}
x(t,{\bf z})&=&x_{0}\cos\,t+y_{0}\sin\,t,\vspace{0.2cm}\\
y(t,{\bf z})&=&y_{0}\cos\,t-x_{0}\sin\,t.
\end{array}
$$

\smallskip

It is clear that the origin of coordinates of $\mathbb{R}^2$ is a global isochronous center for system \eqref{F0_vdp} whose circular periodic solutions starting on ${\bf z}=(x_{0},y_{0})$ $\in\mathbb{R}^{2}$ are pa\-ra\-me\-tri\-zed by the above functions $x(t,{\bf z})$ and $y(t,{\bf z})$. Therefore through every initial condition $(x_{0},y_{0})$ in $\mathbb{R}^2$ passes a $2\pi-$periodic solution of system \eqref{F0_vdp}.

\smallskip

Physically speaking system \eqref{F0_vdp} models a simple harmonic oscillator and its solutions $(x(t,{\bf z}),y(t,{\bf z}))$ describe the wave behaviour of this oscillator. In this direction the problem of perturbation of the global center \eqref{F0_vdp} is equivalent to the problem of perturbation of a simple harmonic oscillator introducing a damping parameter, as external force and considering a potential $V$ by taking small values for $\rho$, $\alpha$, $\lambda$ and $\delta$  as stated in Theorems from 1 to 8.

\smallskip

We note that using the notation of section \ref{averaging} the fundamental matrix $Y(t,{\bf z})$ of system \eqref{F0_vdp} satisfying that $Y(0,{\bf z})$ is the identity of $\mathbb{R}^{2}$ writes
$$
Y(t,{\bf z})=\left(
\begin{array}{cc}
\cos\,t&\sin\,t\vspace{0.2cm}\\
-\sin\,t&\cos\,t
\end{array}
\right).
$$

\smallskip

{From} the averaging theory we are interested in the simple zeros of the function
\begin{equation}\label{average_function}
f({\bf z})=(f_{1}({\bf z}),f_{2}({\bf z}))=\displaystyle\int_{0}^{2\pi}Y^{-1}(t,{\bf z})F_{1}(t,{\bf x}(t,{\bf z},0))dt,
\end{equation}
where $F_{1}(t,{\bf x})$ are the terms of order 1 on $\varepsilon$ of the vector field associated to system \eqref{vdp_new}
$$
(y,-x+\varepsilon^{n_{1}}ry-\varepsilon^{2m+n_{1}}rx^{2}y-\varepsilon^{2m+n_{2}}ax^{3}-\varepsilon^{4m+n_{3}}\ell x^{5}+\varepsilon^{-m+n_{4}}d\cos\,t)^{T}.
$$

\smallskip

Now we prove our results.

\smallskip


\begin{proof}[Proof of Theorem $1$:]
Firstly we take $n_{1}=n_{4}=1$, $m=0$ and consider $n_{2},n_{3}>1$. So the rescaling \eqref{rescaling} becomes $(x,y,\rho,\delta,\alpha,\lambda)=(X,Y,\varepsilon r,\varepsilon d,\varepsilon^{n_{2}} a,$ $\varepsilon^{n_{3}}d)$ and calling again $(X,Y)$ by $(x,y)$ we verify the hypotheses of Theorem 1. Now the vector field of system \eqref{vdp} becomes
$$
(y,-x+\varepsilon ry-\varepsilon rx^{2}y-\varepsilon^{n_{2}}ax^{3}-\varepsilon^{n_{3}}\ell x^{5}+\varepsilon d\cos\,t)^{T}.
$$
Hence as we stated before the terms of order 1 in this vector field are given by $F_{1}(t,{\bf x})=(0,ry-rx^{2}y+d\cos\,t)$. Then since ${\bf z}=(x_{0},y_{0})$ the function $f(x_0,y_0)=(f_{1}(x_0,y_0),f_{2}(x_0,y_0))$ given in \eqref{average_function} turns into the form
$$
\begin{array}{rl}
f_{1}(x_0,y_0)=&\displaystyle\int_{0}^{2\pi}-(\sin\,t)(d\cos\,t+r(y_{0}\cos\,t-x_{0}\sin\,t)-r(y_{0}\cos\,t-\vspace{0.2cm}\\
&\quad\quad x_{0}\sin\,t)(x_{0}\cos\,t+y_{0}\sin\,t)^{2})dt\vspace{0.2cm}\\
=&-\dfrac{1}{4}\pi r x_{0}\left(-4+x_{0}^{2}+y_{0}^{2}\right),\vspace{0.2cm}\\
f_{2}(x_0,y_0)=&\displaystyle\int_{0}^{2\pi}(\cos\,t)(d\cos\,t+r(y_{0}\cos\,t-x_{0}\sin\,t)-r(y_{0}\cos\,t-\vspace{0.2cm}\\
&\quad\quad x_{0}\sin\,t)(x_{0}\cos\,t+y_{0}\sin\,t)^{2})dt\vspace{0.2cm}\\
=&\dfrac{1}{4}\left(4d-r y_{0}\left(-4+x_{0}^{2}+y_{0}^{2}\right)\right).
\end{array}
$$
The zeros of functions $f_{1}$ and $f_{2}$ are a pair of conjugate complex vectors and a real pair $(x_{0}^{0},y_{0}^{0})$ satisfying $x_{0}^{0}=0$ and
$$
y_{0}^{0}=\dfrac{2}{3}\dfrac{(36r)^{1/3}}{\left(9d+\sqrt{81d^{2}-48r^{2}}\right)^{1/3}}+\dfrac{1}{3}\left(\dfrac{6}{r}\right)^{1/3}\left(9d+\sqrt{81d^{2}-48r^{2}}\right)^{1/3}.
$$
Note that $y_{0}^{0}$ is real because $81d^2>48r^2>0$ by assumption. We observe that the expression of $y_{0}^{0}=y_{0}^{0}(r,d)$ does not change when we go back through the rescaling \eqref{rescaling} taking $r=\varepsilon^{-1}\rho$ and $d=\varepsilon^{-1}\delta$ according to the hypotheses of Theorem 1.

\smallskip

Moreover we denote $\Gamma=\left(9d+\sqrt{81d^{2}-48r^{2}}\right)^{1/3}$ and observe that the matrix $M=(m_{ij})=\partial(f_{1},f_{2})/\partial(x_{0},y_{0})$ at $(x_{0}^{0},y_{0}^{0})$ is diagonal and its elements are $m_{11}=-\pi r\Gamma_{+}/36$ and $m_{22}=\pi r\Gamma_{-}$, where
$$
\Gamma_{\pm}=-12\pm\dfrac{6^{1/3}24r^{2/3}}{\Gamma^{2}}\pm\dfrac{6^{2/3}\Gamma^{2}}{r^{2/3}}.
$$

The determinant $\Pi_{0}$ of $M$ written in a power series of $r$ around $r=0$ is
$$
\Pi_{0}=\dfrac{3}{8}\left(  2\pi^{2}d\left(4d\right)^{1/3}\right)r^{2/3}+\dfrac{(\pi r)^{2}}{3}+\mathit{O}(r^{7/3}),  
$$
which is positive for $|r|$ sufficiently small because $d>0$. So the averaging theory described in section \ref{averaging} guarantee the existence of a $2\pi$--periodic solution $(x(t,\varepsilon),y(t,\varepsilon))$ such that $(x(0,\varepsilon),y(0,\varepsilon))$ tends to $(x_{0}^{0},y_{0}^{0})$ when $\varepsilon\to0$.

\smallskip

On the other hand the trace $\Sigma_{0}$ of $M$ at $(x_{0}^{0},y_{0}^{0})$ is
$$
\Sigma_{0}=-2^{2/3}\pi\left(2d\right)^{2/3}r^{1/3}-\dfrac{2\pi r}{3}+\mathit{O}(r^{4/3}).
$$
Thus for $|r|$ sufficiently small the trace $\Sigma_{0}$ of $M$ is negative. So the real part of the eigenvalues of the matrix $M$ are both negative and then the periodic solution $(x(t,\varepsilon),y(t,\varepsilon))$ is stable for each $\varepsilon$ sufficiently small. So Theorem 1 is proved.
\end{proof}

\smallskip


\begin{proof}[Proof of Theorem $2$:]
Now we take $n_{2}=n_{4}=1$ and consider $n_{1},n_{3}>1$. So we are under the hypotheses of Theorem 2. Now, taking $m=0$ and doing the rescaling \eqref{rescaling} the vector field of system \eqref{vdp_new} writes
$$
(y,-x+\varepsilon^{n_{1}} ry-\varepsilon^{n_{1}} rx^{2}y-\varepsilon ax^{3}-\varepsilon^{n_{3}}\ell x^{5}+\varepsilon d\cos\,t)^{T}.
$$
Thus we have $F_{1}(t,{\bf x})=(0,-ax^{3}+d\cos\,t)$ and $f(x_0,y_0)=(f_{1}(x_0,y_0),$ $f_{2}(x_0,y_0))$ becomes
$$
\begin{array}{rcl}
f_{1}(x_0,y_0)&=&\displaystyle\int_{0}^{2\pi}-(\sin\,t)(d\cos\,t-a(x_{0}\cos\,t+y_{0}\sin\,t)^{3}dt\vspace{0.2cm}\\
&=&\dfrac{3}{4}a\pi y_{0}\left(x_{0}^{2}+y_{0}^{2}\right),\vspace{0.2cm}\\
f_{2}(x_0,y_0)&=&\displaystyle\int_{0}^{2\pi}(\cos\,t)(d\cos\,t-a(x_{0}\cos\,t+y_{0}\sin\,t)^{3}dt\vspace{0.2cm}\\
&=&\dfrac{1}{4}\pi\left(     4d-3a x_{0}\left(x_{0}^{2}+y_{0}^{2}\right)    \right).
\end{array}
$$
The real zero $(x_{0}^{0},y_{0}^{0})$ of $f_{1}$ and $f_{2}$ is
$$
(x_{0}^{0},y_{0}^{0})=\left(\left(    \dfrac{4d}{3a}  \right)^{1/3}   ,      0       \right).
$$

Again the values $(x_{0}^{0},y_{0}^{0})$ does not depends on $\varepsilon$ when we go back to the original parameters $\alpha$ and $\delta$ through the rescaling \eqref{rescaling}. Besides the matrix $M=\partial(f_{1},f_{2})/\partial(x_{0},y_{0})$ at $(x_{0}^{0},y_{0}^{0})$ now writes
$$
M=\dfrac{1}{\varepsilon}\left(
\begin{array}{cc}
0&((3/4)\alpha \delta^{2})^{1/3}\pi\vspace{0.2cm}\\
-3((3/4)\alpha \delta^{2})^{1/3}\pi&0
\end{array}
\right),
$$
and then the determinant $\Pi_{0}$ of $M$ is
$$
\Pi_{0}=\dfrac{3\pi^{2}}{2\varepsilon^{2}}\left(  \dfrac{9\alpha^{2}\delta^{4}}{2}   \right)^{1/3}>0.
$$
Therefore the averaging theory implies the existence of a $2\pi$--periodic solution $(x(t,\varepsilon),y(t,\varepsilon))$ such that $(x(0,\varepsilon),y(0,\varepsilon))$ tends to $(x_{0}^{0},y_{0}^{0})$ when $\varepsilon\to0$.
\end{proof}

\smallskip

We observe that since the trace $\Sigma_{0}$ of $M$ is zero and the determinant $\Pi_{0}$ is positive, then the real part of the eigenvalues of $M$ are both zero. So we cannot decide about the stability of the periodic solution of Theorem 2. In what follows we see that the same occurs in Theorem 3.

\smallskip


\begin{proof}[Proof of Theorem $3$:]
We start considering $n_{3}=n_{4}=1$, $n_{1},n_{2}>1$ and $m=0$. Therefore we are under the assumptions of Theorem 3. Doing the rescaling \eqref{rescaling} the vector field of system \eqref{vdp_new} becomes
$$
(y,-x+\varepsilon^{n_{1}} ry-\varepsilon^{n_{1}} rx^{2}y-\varepsilon^{n_{2}} ax^{3}-\varepsilon\ell x^{5}+\varepsilon d\cos\,t)^{T}.
$$
So we conclude that $F_{1}(t,{\bf x})=(0,-\ell x^{5}+d\cos\,t)$ and then we get
$$
\begin{array}{rcl}
f_{1}(x_0,y_0)&=&\displaystyle\int_{0}^{2\pi}-(\sin\,t)(d\cos\,t-\ell(x_{0}\cos\,t+y_{0}\sin\,t)^{5}dt\vspace{0.2cm}\\
&=&\dfrac{5}{8}\ell\pi y_{0}\left(x_{0}^{2}+y_{0}^{2}\right)^{2},\vspace{0.2cm}\\
f_{2}(x_0,y_0)&=&\displaystyle\int_{0}^{2\pi}(\cos\,t)(d\cos\,t-\ell(x_{0}\cos\,t+y_{0}\sin\,t)^{5}dt\vspace{0.2cm}\\
&=&\dfrac{1}{8}\pi\left(     8d-5\ell x_{0}\left(x_{0}^{2}+y_{0}^{2}\right)^{2}    \right).
\end{array}
$$
Now the real zero $(x_{0}^{0},y_{0}^{0})$ of $f_{1}$ and $f_{2}$ is
$$
(x_{0}^{0},y_{0}^{0})=\left(\left(    \dfrac{8d}{5\ell}  \right)^{1/5}   ,      0       \right).
$$

We note that as in Theorems 1 and 2 the root $(x_{0}^{0},y_{0}^{0})$ do not depend on $\varepsilon$ after we perform rescaling \eqref{rescaling}. However using \eqref{rescaling} the matrix $M=\partial(f_{1},f_{2})/\partial(x_{0},y_{0})$ at $(x_{0}^{0},y_{0}^{0})$ is
$$
M=\dfrac{1}{\varepsilon}\left(
\begin{array}{cc}
0&((5/8)\lambda \delta^{4})^{1/5}\pi\vspace{0.2cm}\\
-5((5/8)\lambda \delta^{4})^{1/5}\pi&0
\end{array}
\right),
$$
whose determinant is
$$
\Pi_{0}=\dfrac{5\pi^{2}}{2\varepsilon^{2}}\left(  \dfrac{25\lambda^{2}\delta^{8}}{2}   \right)^{1/5}>0.
$$
Then the averaging method states the existence of a $2\pi$--periodic solution $(x(t,\varepsilon),y(t,\varepsilon))$ such that $(x(0,\varepsilon),y(0,\varepsilon))$ tends to $(x_{0}^{0},y_{0}^{0})$ when $\varepsilon\to0$. The stability of this periodic solution is unknown because the trace $\Sigma_{0}$ of $M$ is zero and the determinant $\Pi_{0}>0$ is positive just as in Theorem 2.
\end{proof}

\smallskip


\begin{proof}[Proof of Theorem $4$:]
In what follows we assume that $n_{1}=m=1$ and $n_{4}=2$. Now the expression of the vector field of system \eqref{vdp_new} is
$$
(y,-x+\varepsilon ry-\varepsilon^{3} rx^{2}y-\varepsilon^{n_{2}+2} ax^{3}-\varepsilon^{n_{3}+4}\ell x^{5}+\varepsilon d\cos\,t)^{T},
$$
and then $F_{1}(t,{\bf x})=(0,ry+d\cos\,t)$. In this case we obtain
$$
\begin{array}{rcl}
f_{1}(x_0,y_0)&=&\displaystyle\int_{0}^{2\pi}-(\sin\,t)(d\cos\,t+r(y_{0}\cos\,t-x_{0}\sin\,t)dt\vspace{0.2cm}\\
&=&\pi r x_{0},\vspace{0.2cm}\\
f_{2}(x_0,y_0)&=&\displaystyle\int_{0}^{2\pi}(\cos\,t)(d\cos\,t+r(y_{0}\cos\,t-x_{0}\sin\,t)dt\vspace{0.2cm}\\
&=&\pi(d+ry_{0}).
\end{array}
$$
It is immediate that the only zero $(x_{0}^{0},y_{0}^{0})$ of $f_{1}$ and $f_{2}$ in this case is $(x_{0}^{0},y_{0}^{0})=(0,d/r)$. So after going back through the rescaling \eqref{rescaling} taking $r=\varepsilon^{-1}\rho$ and $d=\varepsilon^{-2}\delta$ we have
$$
(x_{0}^{0},y_{0}^{0})=\left(       0      ,      -\dfrac{\delta}{\varepsilon \rho}       \right).
$$

In addition the matrix $M=\partial(f_{1},f_{2})/\partial(x_{0},y_{0})$ at $(x_{0}^{0},y_{0}^{0})$ is diagonal and its elements are both given by $\rho\pi/\varepsilon$. Consequently $M$ has determinant $\Pi_{0}=(\rho\pi/\varepsilon)^{2}>0$ and trace $\Sigma_{0}=2\rho\pi/\varepsilon>0$. Then averaging theory described in section \ref{averaging} assure the existence of a $2\pi$--periodic solution $(x(t,\varepsilon),y(t,\varepsilon))$ such that $(x(0,\varepsilon),y(0,\varepsilon))$ tends to $(x_{0}^{0},y_{0}^{0})$ when $\varepsilon\to0$. Moreover since $\rho$ is positive this periodic solution is unstable. Then Theorem 4 is proved.
\end{proof}

\smallskip


\begin{proof}[Proof of Theorem $5$:]
Now we take $m=0$, $n_{2}=n_{3}=n_{4}=1$ and $n_{1}>1$. The expression of the vector field of system \eqref{vdp_new} in this case is
$$
(y,-x+\varepsilon^{n_{1}} ry-\varepsilon^{n_{1}} rx^{2}y-\varepsilon ax^{3}-\varepsilon\ell x^{5}+\varepsilon d\cos\,t)^{T},
$$
and we get $F_{1}(t,{\bf x})=(0,-ax^{3}-\ell x^{5}+d\cos\,t)$. Hence function $f=(f_{1},f_{2})$ writes
$$
\begin{array}{rcl}
f_{1}(x_0,y_0)&=&\displaystyle\int_{0}^{2\pi}-(\sin\,t)(d\cos\,t-a(x_{0}\cos\,t+y_{0}\sin\,t)^{3}-\vspace{0.2cm}\\
&&\quad\quad\ell(x_{0}\cos\,t+y_{0}\sin\,t)^{5}) dt\vspace{0.2cm}\\
&=&\dfrac{1}{8}\pi y_{0}(x_{0}^{2}+y_{0}^{2})(6a+5\ell(x_{0}^{2}+y_{0}^{2})),\vspace{0.2cm}\\
f_{2}(x_0,y_0)&=&\displaystyle\int_{0}^{2\pi}(\cos\,t)(d\cos\,t-a(x_{0}\cos\,t+y_{0}\sin\,t)^{3}-\vspace{0.2cm}\\
&&\quad\quad\ell(x_{0}\cos\,t+y_{0}\sin\,t)^{5}) dt\vspace{0.2cm}\\
&=&\dfrac{1}{8}\pi(8d- x_{0}(x_{0}^{2}+y_{0}^{2})(6a+5\ell(x_{0}^{2}+y_{0}^{2}))).
\end{array}
$$

In order to find the zeros of $f(x_{0},y_{0})$ we compute a
Gr\"{o}bner basis $\{b_{k}(x_{0},$ $y_{0})$, $k=1,2,3\}$ in the variables
$x_{0}$ and $y_{0}$ for the set of polynomials $\{\overline{f}_{1}(x_{0},y_{0}),$ $
\overline{f}_{2}(x_{0},y_{0})\}$, where $\overline{f}_{1,2}= (8/\pi)f_{1,2}$. Then we look for
the zeros of each $b_{k}$, $k=1,2,3$. It is a known fact that the zeros
of a Gr\"{o}bner basis of $\{\overline{f}_{1}(x_{0},y_{0})), \overline{f}_{2}(x_{0},y_{0})\}$ are the zeros of $\overline{f}_{1}$ and $\overline{f}_{2}$, and consequently zeros of $f_{1}$ and $f_{2}$ too.
For more information about a Gr\"{o}bner basis see \cite{AL} and
\cite{L2}.

\smallskip

The Gr\"{o}bner is
$$
\begin{array}{l}
b_{1}(x_{0},y_{0})=dy_{0},\vspace{0.2cm}\\
b_{2}(x_{0},y_{0})=(6ax_{0}^{2}+5\ell x_{0}^{4})y_{0}+(6a+10\ell x_{0}^{2})y_{0}^{3}+5\ell y_{0}^{5},\vspace{0.2cm}\\
b_{3}(x_{0},y_{0})=-8d+6ax_{0}^{3}+5\ell x_{0}^{5}+(6ax_{0}+10\ell x_{0}^{3})y_{0}^{2}+5\ell x_{0}y_{0}^{4}.
\end{array}
$$ 

The only zero of $b_{1}(x_{0},y_{0})$ is $y_{0}^{0}=0$ since $d=\varepsilon^{-1}\delta>0$. So replacing $y_{0}^{0}$ into $b_{3}(x_{0},y_{0})$ and simplifying the new expressions we reduce our problem to find a zero of the polynomial
\begin{equation}\label{pol-c5}
p(x_{0})=\dfrac{1}{5\ell}b_{3}(x_{0},0)=-\dfrac{8d}{r\ell}+\dfrac{6a}{5\ell}x_{0}^{3}+x_{0}^{5}.
\end{equation}

\smallskip

As we know from algebra there is no general formula to provide the roots of a quintic polynomial as polynomial \eqref{pol-c5}. Nevertheless some techniques may become this task easier. An elegant one can be found in \cite{ROOTS}. In this paper it is provide a way to study the number of roots (complex and real, including multiplicity) of a polynomial of any order only performing some calculations on the coefficients of the considered polynomial. We will present a brief summary of the algorithmic of this method in the Appendix \ref{appendix} at the end of this paper.

\smallskip

By the Fundamental Theorem of Algebra polynomial \eqref{pol-c5} has 5 roots taking into account their multiplicities and once we want to apply the averaging theory we are only interested in the real simple roots of polynomial \eqref{pol-c5}. Indeed, applying the method due to \cite{ROOTS} and summarized in Appendix \ref{appendix} we obtain
$$
\begin{array}{rl}
D_{4}=&\dfrac{384ad^{2}}{\ell^{3}},\vspace{0.2cm}\\
D_{5}=&\dfrac{53747712}{78125}\dfrac{a^5 d^2}{\ell^7}+20480\dfrac{d^4}{\ell^4}.
\end{array}
$$
Therefore as described in Appendix \ref{appendix}, since $\alpha\lambda=\varepsilon^{2}a\ell<0$ by hypothesis we have that $D_{4}$ is negative. So polynomial \eqref{pol-c5} has a unique simple real root if $D_{5}=D=0$ according to statement $(5)$ of Appendix \ref{appendix}. On the other hand statement $(3)$ says that polynomial \eqref{pol-c5} has three simple real roots if $D_{5}=D<0$. Moreover $D$ and $D_{4}$ does not change when we replace $(a,\ell,d)$ by $(\varepsilon^{-1}\alpha,\varepsilon^{-1}\lambda,\varepsilon^{-1}\delta)$ using rescaling \eqref{rescaling}.

\smallskip

Now using the property on zeros of the Gr\"{o}bner basis cited previously the function $f(x_{0},y_{0})$ has a zero $(x_{0}^{0},0)$ if $D=0$ and three zeros $(x_{0}^{i},0)$ if $D<0$, with $i=1,2,3$. Besides the matrix $M=\partial(f_{1},f_{2})/\partial(x_{0},y_{0})$ evaluated at $y_{0}^{0}=0$ after rescaling \eqref{rescaling} is
$$
M=\dfrac{1}{\varepsilon}\left(
\begin{array}{cc}
0&\dfrac{\pi}{8}(6\alpha x_{0}^{2}+5\lambda x_{0}^{4})\vspace{0.2cm}\\
-\dfrac{\pi}{8}(18\alpha x_{0}^{2}+25\lambda x_{0}^{4})&0
\end{array}
\right),
$$
whose determinant is
$$
\Pi(x_{0})=\dfrac{1}{\varepsilon^{2}}\left(\dfrac{27}{16}(\alpha\pi x_{0}^{2})^{2}+\dfrac{15}{4}\alpha\lambda(\pi x_{0}^{3})^{2}+\dfrac{125}{64}(\lambda\pi x_{0}^{4})^{2}\right).
$$
Thus by computing a Gr\"{o}bner basis for polynomials $\{\Pi(x_{0}),
f_{1}(x_{0},0),f_{2}(x_{0},$ $0)\}$ in the variables $\alpha$, $\delta$, $\lambda$ and $x_{0}$ we get the polynomials $\overline{b}_{1}(\alpha,\delta,\lambda,x_{0})=\delta$ and $\overline{b}_{2}(\alpha,\delta,\lambda,x_{0})=6\alpha+5\lambda x_{0}^{2}$. So $\overline{b}_{1}(\alpha,\delta,\lambda,x_{0})>0$ because $\delta>0$ and it means that we cannot have $f_{1}$, $f_{2}$ and $\Pi$ equal to zero simultaneously. Consequently $\Pi(x_{0}^{i})\neq0$ for all $i=0,1,2,3$. Then using the averaging theory we conclude the existence of a $2\pi$--periodic solution if $D=0$ and three $2\pi$--periodic solution if $D<0$.

\smallskip

Now we exhibit values of $\alpha$, $\lambda$ and $\delta$ for which system \eqref{vdp} has one or three $2\pi$--periodic solutions. Despite, if consider $a=25/18$, $\ell=-1$ and $d=5/12$ we obtain $a\ell=-25/18<0$ and $D=0$. On the other hand, taking $d=18a/31$ and $\ell=-42a/155$ we get $a\ell=-42a^{2}/155<0$ and $D=-1425339825408/823543<0$. So Theorem 5 is proved.
\end{proof}

\smallskip

We note that by assumption in Theorem 5 we have $a\ell<0$. So considering the notation of Appendix \ref{appendix} the coefficients $D_{2}=-6a/5\ell$ and $D_{3}=384ad^2/\ell^3$ corresponding to polynomial \eqref{pol-c5} are positive and $D_{4}$ is negative. Thus the only possible configurations of roots of the polynomial \eqref{pol-c5} are those one listed in statements $(3)$ and $(5)$ of the Appendix. This fact means that we cannot have 2 or 5 periodic solutions when we consider $m=0$, $n_{2}=n_{3}=n_{4}=1$ and $n_{1}>1$.

\smallskip

%
%


\begin{proof}[Proof of Theorem $6$:]
We start fixing the values $m=0$, $n_{1}=n_{2}=n_{4}=1$ and $n_{3}>1$. Now the vector field of system \eqref{vdp_new} writes
$$
(y,-x+\varepsilon ry-\varepsilon rx^{2}y-\varepsilon ax^{3}-\varepsilon^{n_{3}}\ell x^{5}+\varepsilon d\cos\,t)^{T},
$$
and then $F_{1}(t,{\bf x})=(0,ry- rx^{2}y-ax^{3}+d\cos\,t)$. With this expression of $F_{1}$ function \eqref{average_function} is
$$
\begin{array}{rcl}
f_{1}(x_0,y_0)&=&\displaystyle\int_{0}^{2\pi}-(\sin\,t)(d\cos\,t+r(y_{0}\cos\,t-x_{0}\sin\,t)\vspace{0.2cm}\\
&&\quad\quad-r(y_{0}\cos\,t-x_{0}\sin\,t)(x_{0}\cos\,t+y_{0}\sin\,t)^{2} \vspace{0.2cm}\\
&&\quad\quad-a(x_{0}\cos\,t+y_{0}\sin\,t)^{3})dt\vspace{0.2cm}\\
&=&-\dfrac{1}{4}\pi(rx_{0}(-4+x_{0}^{2}+y_{0}^{2})-3ay_{0}(x_{0}^{2}+y_{0}^{2})),\vspace{0.2cm}\\
f_{2}(x_0,y_0)&=&\displaystyle\int_{0}^{2\pi}(\cos\,t)(d\cos\,t+r(y_{0}\cos\,t-x_{0}\sin\,t)\vspace{0.2cm}\\
&&\quad\quad-r(y_{0}\cos\,t-x_{0}\sin\,t)(x_{0}\cos\,t+y_{0}\sin\,t)^{2} \vspace{0.2cm}\\
&&\quad\quad-a(x_{0}\cos\,t+y_{0}\sin\,t)^{3})dt\vspace{0.2cm}\\
&=&\dfrac{1}{4}\pi(4d-ry_{0}(-4+x_{0}^{2}+y_{0}^{2})-3ax_{0}(x_{0}^{2}+y_{0}^{2})).
\end{array}
$$

In order to find zeros $(x_{0}^{*},y_{0}^{*})$ of $f(x_{0},y_{0})$ as before we compute a
Gr\"{o}bner basis $\{b_{k}(x_{0},y_{0})$, $k=1,\ldots,14\}$ in the variables
$x_{0}$ and $y_{0}$ now for the set of polynomials $\{\overline{f}_{1}(x_{0},y_{0}),
\overline{f}_{2}(x_{0},y_{0})\}$ with $\overline{f}_{1,2}= \mp(4/\pi)f_{1,2}$. We will look for
zeros of two elements of the Gr\"{o}bner basis for the polynomials $\{\overline{f}_{1}(x_{0},y_{0}),
\overline{f}_{2}(x_{0},y_{0})\}$ in the variables $x_{0}$ and $y_{0}$. These polynomials are
$$
\begin{array}{rl}
b_{1}(x_{0},y_{0})=&144 a^2 d r^3-4 d^3 r^3+(108 a^2 d^2 r^2+144 a^2 r^4-4 d^2 r^4)y_{0}+\vspace{0.2cm}\\
&144 a^2 d r^3y_{0}^{2}+(81 a^4 d^2+18 a^2 d^2 r^2+d^2 r^4)y_{0}^{3},\vspace{0.2cm}\\
b_{2}(x_{0},y_{0})=&216 a^3 d r^2+(-324 a^4 r^2-9 a^2 d^2 r^2+36 a^2 r^4-d^2 r^4)x_{0}+\vspace{0.2cm}\\
&(27 a^3 d^2 r + 216 a^3 r^3 + 3 a d^2 r^3)y_{0}+(243 a^5 d + 54 a^3 d r^2 +\vspace{0.2cm}\\
&3 a d r^4)y_{0}^{2}.
\end{array}
$$ 

We must observe that $b_{1}(x_{0},y_{0})$ depends only on $y_{0}$. Then for each zero $y_{0}^{*}$ of $b_{1}(y_{0})$ the second polynomial $b_{2}(x_{0},y_{0})$ provides a zero $x_{0}^{*}$ associated to $y_{0}^{*}$ because the coefficient of $x_{0}$ in $b_{2}(x_{0},y_{0})$ is not zero by hypothesis. Now we will look for zeros of $b_{1}(y_{0})$. Indeed the discriminant $\Delta$ of the cubic polynomial $b_{1}(y_{0})$ writes
$$
\begin{array}{rl}
\Delta=&-16 d^2 r^6\left(324 a^4+d^2 r^2+9 a^2\left(d^2-4 r^2\right)\right)^{2}\left(2187 a^4 d^4+27 d^4 r^4-\right.\vspace{0.2cm}\\
&\left.16 d^2 r^6+18 a^2\left(27 d^4 r^2-72 d^2 r^4+32 r^6\right)\right)\vspace{0.2cm}\\
=&-16 d^2 r^6\Delta_{1}^{2}\Delta_{2}.
\end{array}
$$
So if $\Delta_{1}$ or $\Delta_{2}$ is zero then $\Delta$ is zero and consequently $b_{1}(y_{0})$ has a simple real root $y_{0}^{0}$ because this polynomial has no root with multiplicity 3. Actually $b_{1}'''(y_{0})=6 (81 a^4 d^2 + 18 a^2 d^2 r^2 + d^2 r^4)>0$. The same occurs if $\Delta_{2}>0$ since this condition implies $\Delta<0$. On the other hand, if $\Delta_{2}<0$ then $\Delta>0$ and consequently polynomial $b_{1}(y_{0})$ has three simple real roots $y_{0}^{i}$, $i=1,2,3$. Additionally, replacing each value $y_{0}^{i}$ into $b_{2}(x_{0},y_{0})$ we obtain the respective values $x_{0}^{i}$, for $=0,1,2,3$.  We note that coming back through the rescaling \eqref{rescaling} the signs of $\Delta$, $\Delta_{1}$ and $\Delta_{2}$ does not change because each monomial composing $\Delta$ has the same degree.

\smallskip

Now we will verify the condition $M=\det((\partial f/\partial {\bf z})(x_{0}^{i},y_{0}^{i}))\neq0$ for $i=0,1,2,3$. In fact $M$ is
$$
\dfrac{-\pi}{4\varepsilon}\left(
\begin{array}{cc}
 -6\alpha x_{0}y_{0}+\rho\left(  -4+3x_{0}^{2}+y_{0}^{2}  \right)& 2\rho x_{0}y_{0}-3\alpha\left(  x_{0}^{2}+3y_{0}^{2}  \right)\vspace{0.2cm}\\
2\rho x_{0}y_{0}+3\alpha\left(  3x_{0}^{2}+y_{0}^{2}   \right)& 6\alpha x_{0}y_{0}+\rho\left(  -4+x_{0}^{2}+3y_{0}^{2}   \right)
\end{array}
\right),
$$
whose determinant $\Pi$ now is
$$
\Pi(x_{0},y_{0})=\dfrac{\pi^{2}}{16\varepsilon^{2}}\left(    27\alpha^{2}\left(x_{0}^{2}+y_{0}^{2}\right)^{2}+\rho^{2}\left(  -4+x_{0}^{2}+y_{0}^{2}   \right)\left(  -4+3x_{0}^{2}+3y_{0}^{2}   \right)    \right).
$$
However for each $i=0,1,2,3$ the determinant $\Pi(x_{0}^{i},y_{0}^{i})\neq0$ by hypothesis. Therefore using averaging theory described in section \ref{averaging} system \eqref{vdp} has a $2\pi$--periodic solution if $\Delta_{1}\Delta_{2}=0$ or $\Delta_{2}>0$ and three $2\pi$--periodic solution if $\Delta_{2}<0$.

\smallskip

Now we present concrete values of $\alpha$, $\rho$ and $\delta$ for which we have one or three periodic solutions. Indeed taking $a=r=1$ and $d=6$ we have $\Delta_{2}=3452544>0$. Meantime if we consider $d=6a$ and
$$
r=-\dfrac{9a}{\sqrt{-9+4\sqrt{6}}},
$$
we obtain $\Delta_{2}=-34012224 a^{8}/(9-4\sqrt{6})^{2}<0$. We also observe that the values $(r,a,d)=((555/4)\sqrt{\dfrac{154073}{9622}},185,22)$ and $(r,a,d)=(585\sqrt{3},195,$ $585\sqrt{2})$ make $\Delta_{1}$ and $\Delta_{2}$ equal to zero, respectively.
This ends the proof of Theorem 6.
\end{proof}

\smallskip


\begin{proof}[Proof of Theorem $7$:]
For this case we assume $m=0$ and $n_{i}=1$, $i=1,2,3,4$. Then the vector field of system \eqref{vdp_new} writes
$$
(y,-x+\varepsilon ry-\varepsilon rx^{2}y-\varepsilon ax^{3}-\varepsilon\ell x^{5}+\varepsilon d\cos\,t)^{T},
$$
and we have $F_{1}(t,{\bf x})=(0,ry-rx^{2}y-ax^{3}-\ell x^{5}+d\cos\,t)$. So it follows that
$$
\begin{array}{rl}
f_{1}(x_0,y_0)=&\displaystyle\int_{0}^{2\pi}-(\sin\,t)(d\cos\,t+r(y_{0}\cos\,t-x_{0}\sin\,t)-r(y_{0}\cos\,t-\vspace{0.2cm}\\
&\quad\quad x_{0}\sin\,t)(x_{0}\cos\,t+y_{0}\sin\,t)^{2}-a(x_{0}\cos\,t+y_{0}\sin\,t)^{3}-\vspace{0.2cm}\\
&\quad\quad-\ell(x_{0}\cos\,t+y_{0}\sin\,t)^{5}) dt\vspace{0.2cm}\\
=&-\dfrac{1}{8}\pi (2rx_{0}(-4+x_{0}^{2}+y_{0}^{2})-y_{0}(x_{0}^{2}+y_{0}^{2})(6a+5\ell(x_{0}^{2}+y_{0}^{2})),\vspace{0.2cm}\\
f_{2}(x_0,y_0)=&\displaystyle\int_{0}^{2\pi}(\cos\,t)(d\cos\,t+r(y_{0}\cos\,t-x_{0}\sin\,t)-r(y_{0}\cos\,t-\vspace{0.2cm}\\
&\quad\quad x_{0}\sin\,t)(x_{0}\cos\,t+y_{0}\sin\,t)^{2}-a(x_{0}\cos\,t+y_{0}\sin\,t)^{3}-\vspace{0.2cm}\\
&\quad\quad-\ell(x_{0}\cos\,t+y_{0}\sin\,t)^{5}) dt\vspace{0.2cm}\\
=&\dfrac{1}{8}\pi(8d+8ry_{0}-(x_{0}^{2}+y_{0}^{2})(6ax_{0}+2ry_{0}+5\ell x_{0}(x_{0}^{2}+y_{0}^{2}))).
\end{array}
$$

We will find zeros of $f(x_0,y_0)$ through the zeros of a Gr\"{o}bner basis $\{b_{k}(x_{0},y_{0})\}$ of it. As before instead of $f(x_0,y_0)$ the Gr\"{o}bner basis will be related with the functions $\overline{f}_{1,2}=\mp(8/\pi)f_{1,2}$.

\smallskip

One of the elements of the Gr\"{o}bner basis in the variables $x_{0}$ and $y_{0}$ is
$$
\begin{array}{rl}
b_{1}(x_{0},y_{0})=&144 a^2 d r^5 - 4 d^3 r^5 + 960 a d \ell r^5 + 1600 d \ell^2 r^5+(108 a^2 d^2 r^4 +\vspace{0.2cm}\\
& 960 a d^2 \ell r^4 + 2000 d^2 \ell^2 r^4 + 144 a^2 r^6 -4 d^2 r^6 + 960 a \ell r^6 +\vspace{0.2cm}\\
&1600 \ell^2 r^6)y_{0}+(120 a d^3 \ell r^3 + 500 d^3 \ell^2 r^3 + 144 a^2 d r^5 + 1200\vspace{0.2cm}\\
&a d \ell r^5 + 2400 d \ell^2 r^5)y_{0}^{2}+81 a^4 d^2 r^2 + 540 a^3 d^2 \ell r^2 + 900 a^2 d^2 \ell^2 \vspace{0.2cm}\\
&r^2 + 18 a^2 d^2 r^4 + 300 a d^2 \ell r^4 + 900 d^2 \ell^2 r^4 + d^2 r^6)y_{0}^{3}+(-450  \vspace{0.2cm}\\ 
&a^2 d^3 \ell^2 r- 1500 a d^3 \ell^3 r + 50 d^3 \ell^2 r^3)y_{0}^4 + 625 d^4 \ell^4 y_{0}^5.
\end{array}
$$ 

We note that again $b_{1}(x_{0},y_{0})$ depends only on $y_{0}$. Moreover, the Gr\"{o}bner basis has other element $b_{2}(x_{0},y_{0})$ which is linear on $x_{0}$ and the coefficients depend on $y_{0}$, where the coefficient $C$ of $x_{0}$ is not zero by hypotheses. Its expression is too large and we will omit it here. It can be easily obtained by an algebraic manipulator. So for each zero $y_{0}^{*}$ of $b_{1}(y_{0})$ the second polynomial $b_{2}(x_{0},y_{0})$ provides a zero $x_{0}^{*}$ associated to $y_{0}^{*}$.

\smallskip

Now we will look for zeros of $b_{1}(y_{0})$. Indeed using again Appendix \ref{appendix} the conditions on $a$, $r$, $d$ and $\ell$ that provide three zeros, $(x_{0}^{i},y_{0}^{i})$, where $i=1,2,3$ for $b_{1}(y_{0})$ are $D_{5}=-D_{5}^{1}D_{5}^{2}D<0$, where $D$ is given in Theorem 7 and $D_{5}^{1}$ and $D_{5}^{2}$ are the values
$$
\begin{array}{ll}
D_{5}^{1}=&\dfrac{110075314176a^{8}r^{20}}{d^{20}(-9a^{2}+r^{2})^{28}},\vspace{0.2cm}\\
D_{5}^{2}=&(59049 a^{10} d^4+16 r^{14}+36 a^2 r^{10}(11 d^2 - 4 r^2)-81 a^4 r^6(d^4 - 12 d^2 r^2 +\vspace{0.2cm}\\
&16 r^4)-6561 a^8(3 d^4 r^2 - 4 d^2 r^4)+729 a^6(3 d^4 r^4 - 28 d^2 r^6 + 16 r^8))^{2}.
\end{array}
$$

By hypothesis $ar=\varepsilon^{-2}\alpha\rho\neq0$ and $D_{5}^{2}\neq0$, then $D_{5}^{1}$ and $D_{5}^{2}$ are positive. Therefore, in order to have $D_{5}$ negative we must have $D$ positive. On the other hand, if $D_{5}$ is positive then according to statement $(2)$ of Appendix \ref{appendix} we have a unique zero $(x_{0}^{0},y_{0}^{0})$ for $b_{1}(y_{0})$ because the value 
$$
D_{2}=-\dfrac{1296 a^{2}r^4 (-3 a^2 + r^2)^{2}}{d^2 (-9 a^2 + r^{2})^{4}},
$$
is negative if $(r^{2}-3a^{2})(r^{2}-9a^{2})\neq0$. In addition by considering $a=\varepsilon^{-1}\alpha$, $r=\varepsilon^{-1}\rho$ and $d=\varepsilon^{-1}\delta$ we get
$$
M=\dfrac{-\pi}{8\varepsilon}\left(
\begin{array}{cc}
2\rho(-4+3x_{0}^{2}+y_{0}^{2})-\Gamma_{1}&\Gamma_{2}-\Gamma_{4}-5\lambda\Gamma_{3}(x_{0}^{2}+5y_{0}^{2})\\
\Gamma_{2}+\Gamma_{4}+5\lambda\Gamma_{3}(5x_{0}^{2}+y_{0}^{2})&2\rho(-4+x_{0}^{2}+3y_{0}^{2})-\Gamma_{1}
\end{array}
\right)
$$
where $\Gamma_{1}=2 x_{0} y_{0} (3 \alpha + 5 \ell (x_{0}^2 + y_{0}^2))$, $\Gamma_{2}=-4\rho x_{0}y_{0}$, $\Gamma_{3}=x_{0}^{2}+y_{0}^{2}$ and $\Gamma_{4}=6\alpha(x_{0}^{2}+	3y_{0}^{2})$. Moreover the determinant $\Pi$ of $M$ writes
$$
\begin{array}{rl}
\Pi(x_{0},y_{0})=&\dfrac{1}{\varepsilon^{2}}\left(4\rho^{2}\left(-4+x_{0}^{2}+y_{0}^{2}\right)(-4+3x_{0}^{2}+3y_{0}^{2})+(x_{0}^{2}+y_{0}^{2})^{2}(6\alpha+\right.\vspace{0.2cm}\\
&\left.\rho\lambda(x_{0}^{2}+y_{0}^{2})(18\alpha+25\lambda x_{0}^{2}+3y_{0}^{2}))\right).
\end{array}
$$
However, by hypothesis $\Pi(x_{0}^{i},y_{0}^{i})\neq0$ for $i=0,1,2,3$ and hence we conclude the first part of Theorem 7.

\smallskip

Now we provide values for $\alpha$, $\rho$, $\lambda$ and $\delta$ for which we have one or three $2\pi$--periodic solutions. First we observe that taking $(a,r,d)=(1,2,8/3)$ we obtain
$$
D=-\dfrac{2770035802112}{6561}<0.
$$

However if we consider
$$
(a,r,d,\ell)=(1/\sqrt{3},1,1,-1/(5\sqrt{3})),
$$
then $D=2752$. So Theorem 7 is proved.
\end{proof}

\smallskip

Now we will prove that under the conditions $m=0$ and $n_{i}=1$, $i=1,2,3,4$ and the assumptions of Theorem 7 we cannot have 2 or 5 periodic solutions by using the averaging theory and the rescaling \eqref{rescaling}. Indeed we observe that the only possibility to have 5 periodic solutions that the statement $(1)$ of Appendix \ref{appendix} holds. However this condition does not hold since $D_{2}$ is negative for any $a$, $r$ and $d$ satisfying the hypotheses of Theorem 7.

\smallskip

In addition we note that the only condition that provides two periodic solutions according with Appendix \ref{appendix} is statement $(7)$, and it needs $D_{5}=D_{4}=E_{2}=0$. Nevertheless we start considering this condition and we will find a Gr\"{o}bner basis for the set of polynomials $\{D_{4},D_{5},E_{2}\}$ in the variables $a$, $d$ and $r$. First we will consider the factor $D_{5}^{2}$ of $D_{5}$ because $D_{5}^{1}\neq 0$ by hypothesis. In this case we obtain 55 polynomials in the Gr\"{o}bner basis, where the first one is $g(d,r)=r^{62}\overline{g}(d,r)$ and
$$
\begin{array}{rl}
\overline{g}(d,r)=&91125 d^{12} - 438800 d^{10} r^2 + 919360 d^8 r^4 - 1038464 d^6 r^6 + 
 \vspace{0.2cm}\\
 & 660992 d^4 r^8 -225280 d^2 r^{10} + 32768 r^{12}.
\end{array}
$$


\smallskip

Since $r\neq0$ we must find a zero of $\overline{g}(d,r)$. So solving $\overline{g}(d,r)=0$ in the variable $r$ we obtain 6 pairs of values
$$
r_{\pm}^{i}=\pm d\sqrt{q_{0}^{i}},
$$
where $i=1,\ldots,6$ and for each $i$ the value $q_{0}^{i}$ is the $i{\it th}-$root of the polynomial
$$
\begin{array}{rl}
q(x)=& 91125 - 438800 x + 919360 x^2 - 1038464 x^3 + 660992 x^4  \vspace{0.2cm}\\
 &- 
  225280 x^5 +32768 x^6.
\end{array}
$$
 
This polynomial has only complex roots. Indeed, if we study the function $q'(x)$ which is given by
$$
\begin{array}{rl}
q'(x)=&-438800 + 1838720 x - 3115392 x^2 + 2643968 x^3 - 1126400 x^4\vspace{0.2cm}\\
& + 
 196608 x^5,
\end{array}
$$
we can apply the method described in Appendix \ref{appendix} to show that the respective values $D_{2}$ and $D_{5}$ related to $q'(x)$ are negative and positive, respectively. Therefore $q'(x)$ has only one real root $q'_{0}$ whose approximate value is $q'_{0}=106703/100000$. Additionally the approximate value of $q(x)$ evaluated in $q'_{0}$ is $q(q'_{0})=16558687/10000>0$. So if we observe that the coefficient of $x^{6}$ is positive it follows that the minimum value of $q(x)$ is positive. Then each $r_{\pm}^{i}$, $i=1,\ldots,6$ is a complex value and the correspondent Gr\"{o}bner basis has no zeros. Then we cannot have a common real zero of $D_{4}$, $D_{5}^{2}$ and $E_{2}$ and so statement $(7)$ in Appendix \ref{appendix} does not hold. Consequently it is not possible to obtain 2 periodic solutions. The proof considering the factor $D$ of $D_{5}$ instead of $D_{5}^{2}$ leads to the same polynomial $g(d,r)$ and then we have the same conclusion.

\smallskip


\begin{proof}[Proof of Theorem $8$:]
In what follows we take $m=0$, $n_{1}=n_{3}=n_{4}=1$ and $n_{2}>1$. So the vector field of system \eqref{vdp_new} becomes
$$
(y,-x+\varepsilon ry-\varepsilon rx^{2}y-\varepsilon^{n_{2}}ax^{3}-\varepsilon\ell x^{5}+\varepsilon d\cos\,t)^{T}.
$$
and we obtain $F_{1}(t,{\bf x})=(0,ry-rx^{2}y-\ell x^{5}+d\cos\,t)$. So $f(x_{0},y_{0})$ writes
$$
\begin{array}{rl}
f_{1}(x_0,y_0)=&\displaystyle\int_{0}^{2\pi}-(\sin\,t)(d\cos\,t+r(y_{0}\cos\,t-x_{0}\sin\,t)-r(y_{0}\cos\,t-\vspace{0.2cm}\\
&\quad\quad x_{0}\sin\,t)(x_{0}\cos\,t+y_{0}\sin\,t)^{2}-\ell(x_{0}\cos\,t+y_{0}\sin\,t)^{5}) dt\vspace{0.2cm}\\
=&-\dfrac{1}{8}\pi (2rx_{0}(-4+x_{0}^{2}+y_{0}^{2})-5\ell y_{0}(x_{0}^{2}+y_{0}^{2})^{2}),\vspace{0.2cm}\\
f_{2}(x_0,y_0)=&\displaystyle\int_{0}^{2\pi}(\cos\,t)(d\cos\,t+r(y_{0}\cos\,t-x_{0}\sin\,t)-r(y_{0}\cos\,t-\vspace{0.2cm}\\
&\quad\quad x_{0}\sin\,t)(x_{0}\cos\,t+y_{0}\sin\,t)^{2}-\ell(x_{0}\cos\,t+y_{0}\sin\,t)^{5}) dt\vspace{0.2cm}\\
=&\dfrac{1}{8}\pi(8d-2ry_{0}(-4+x_{0}^{2}+y_{0}^{2})-5\ell x_{0}(x_{0}^{2}+y_{0}^{2})^{2}).
\end{array}
$$

Again we will find zeros of $f(x_0,y_0)$ through the roots of a Gr\"{o}bner basis $\{b_{k}(x_{0},y_{0})\}$ of $f(x_{0},y_{0})$. We will find a Gr\"{o}bner basis of the functions $\overline{f}_{1,2}$ $=\mp(8/\pi)f_{1,2}$.

\smallskip

An element of the Gr\"{o}bner basis in the variables $x_{0}$ and $y_{0}$ is
$$
\begin{array}{rl}
b_{1}(x_{0},y_{0})=&-4 d^3 r^5 + 1600 d L^2 r^5 +(2000 d^2 \ell^2 r^4 - 4 d^2 r^6+1600 L^2 R^6)y_{0}    \vspace{0.2cm}\\
&  +(500 d^3 \ell^2 r^3 + 2400 d \ell^2 r^5)y_{0}^{2}+(900 d^2 \ell^2 r^4 + d^2 r^6)y_{0}^{3} +     \vspace{0.2cm}\\
&50 d^3\ell^2 r^3y_{0}^4 + 625 d^4 \ell^4 y_{0}^5.
\end{array}
$$ 

As before $b_{1}(x_{0},y_{0})$ depends only on $y_{0}$ and another element $b_{2}(x_{0},y_{0})$ of the Gr\"{o}bner basis is linear on $x_{0}$ with coefficients depending on $y_{0}$,  where the coefficient $C$ of $x_{0}$ in $b_{2}(x_{0},y_{0})$ is not zero by hypotheses. We will not present the expression of $b_{2}(x_{0},y_{0})$ in order to avoid large expressions. Then for each zero $y_{0}^{*}$ of $b_{1}(y_{0})$ we have a second zero $x_{0}^{*}$ through $b_{2}(x_{0},y_{0})$ related to $y_{0}^{*}$.

\smallskip

In order to apply the method described in \cite{ROOTS} and Appendix \ref{appendix} for $b_{1}(y_{0})$ we will perform the translation $y_{0}=\varphi-(2r^{3}/d(125\ell)^{2})$. With this translation we obtain a new polynomial $b_{1}^{*}$ into the form
$$
b_{1}^{*}(\varphi)=\varphi^{5}+\dfrac{4500 \ell^2 r^4 - 3 r^6}{3125 d^2 \ell^4}\varphi^{3}+\dfrac{2k_{1}}{\overline{k_{1}}d^{3}\ell^{6}}\varphi^{2}+\dfrac{4r^{4}k_{2}}{\overline{k_{2}}d^{4}\ell^{8}}\varphi-\dfrac{4r^{5}k_{3}}{\overline{k_{3}}d^{5}\ell^{10}},
$$
where
$$
\begin{array}{rl}
k_{1}=&156250 d^2 \ell^4 r^3 + 750000 \ell^4 r^5 - 13500 \ell^2 r^7 + r^9,\vspace{0.2cm}\\
\overline{k_{1}}=&390625,\vspace{0.2cm}\\
k_{2}=&39062500 d^2 \ell^6 - 390625 \ell^4 (d^2 - 80 \ell^2) r^2 - 1500000 \ell^4 r^4 + 
 13500 \ell^2 r^6 +\vspace{0.2cm}\\
& 3 r^8,\vspace{0.2cm}\\
\overline{k_{2}}=&48828125,\vspace{0.2cm}\\
k_{3}=&48828125 d^4 \ell^6 - 781250 d^2 (25000 \ell^8 - 500 \ell^6 r^2 + 3 \ell^4 r^4) + 
 2 r^4 \vspace{0.2cm}\\
 &(156250000 \ell^6 - 3750000 \ell^4 r^2 + 22500 \ell^2 r^4 + 9 r^6),\vspace{0.2cm}\\
\overline{k_{3}}=&30517578125 .
\end{array}
$$

We observe that the translation performed do not change the number or kind of the zeros of the original polynomial $b_{1}$. Now we will apply the method of Appendix \ref{appendix} for $b_{1}^{*}$. Indeed we have
$$
\begin{array}{rl}
D_{2}=&N_{2},\vspace{0.2cm}\\
D_{3}=&-\dfrac{48r^{6}}{9765625 d^6 \ell^{10}}N_{3},\vspace{0.2cm}\\
D_{4}=&-\dfrac{16r^{12}}{3814697265625 d^{12} \ell^{18}}N_{4},\vspace{0.2cm}\\
D_{5}=&\dfrac{256r^{20}}{1490116119384765625 d^{20} \ell^{26}}N_{5}M_{5}.
\end{array}
$$
Therefore if $D_{5}$ is negative then $b_{1}^{*}$ and consequently $b_{1}$ has exactly three zeros and hence function $f=(f_{1},f_{2})$ has exactly three zeros $(x_{0}^{i},y_{0}^{i})$, $i=1,2,3$. On the other hand if $D_{5}$ is positive and one of the values $D_{2}$, $D_{3}$ or $D_{4}$ is non--positive, then $f$ has exactly one zero $(x_{0}^{0},y_{0}^{0})$. Moreover using the rescaling \eqref{rescaling} the matrix $M$ is the same than the one of the proof of Theorem 7 taking $a=0$. In addition the determinant $\Pi$ of $M$ writes
$$
\Pi(x_{0},y_{0})=\dfrac{1}{\varepsilon^{2}}\left(125\lambda^{2}\left(x_{0}^{2}+y_{0}^{4}\right)+4\rho^{2}(-4+x_{0}^{2}+y_{0}^{2})(-4+3x_{0}^{2}+3y_{0}^{2})\right).
$$
By hypothesis $\Pi(x_{0}^{i},y_{0}^{i})\neq0$ for $i=0,1,2,3$, then we have the first part of Theorem 8 proved.

\smallskip

Now we exhibit values of $\rho$, $\lambda$ and $\delta$ for which we have either one, or three periodic solutions. Firstly we consider
$$
(r,\ell,d)=\left(  \sqrt{31}   ,   1   ,   \dfrac{1}{125}\sqrt{1547875 - 226851\sqrt{35} }  \right).
$$

With these values we obtain
$$
\begin{array}{l}
N_{3}=\dfrac{1}{625}(206837701817900 - 10419995302263\sqrt{35}),\vspace{0.2cm}\\
N_{5}=\dfrac{2(-649491478051728715458244 + 109639832554794027558525\sqrt{35})}{48828125},
\end{array}
$$
which are positive. On the other hand, considering $(r,\ell,d)=(\sqrt{31}/4,1,$ $1)$ we get $N_{5}=-258261575$ $/1048576<0$. This ends the proof of Theorem 8.
\end{proof}

\smallskip

Now we will show that if $m=0$, $n_{1}=n_{3}=n_{4}=1$ and $n_{2}>1$, then we cannot have 2 periodic solutions by using averaging theory described in section \ref{averaging} and rescaling \eqref{rescaling}. Indeed in order that the function $b_{1}^{*}(\varphi)$, given in the proof of Theorem 8, has two periodic solutions we need statement $(7)$ of Appendix \ref{appendix}. We start considering the conditions $N_{5}=N_{4}=E_{2}=0$ and we will see that these equalities implies that $N_{3}$ cannot be negative. We claim that $N_{5}$ and $M_{5}$ are factors of $D_{5}$ as we saw in the last proof. Also we will not present the expression of $E_{2}$ in order to avoid its large expression. It can be obtained by an algebraic manipulator as Mathematica. Additionally we note that $N_{4}=N_{5}=0$ implies $D_{4}=D_{5}=0$, which are necessary conditions for having exactly 2 real roots.

\smallskip

We start replacing the condition $N_{4}=N_{5}=E_{2}=0$ for other ones easier when $r>0$. It is easy to see that if $r=0$ we have $N_{3}=0$ and consequently we cannot get 2 roots for $b_{1}^{*}$. Now we obtain a Gr\"{o}bner basis for the set of polynomials $\{N_{4},N_{5},E_{2}\}$ in the variables $d$, $r$ and $\ell$. This basis has twelve elements, where two of them are
$$
\begin{array}{rl}
g_{1}(r,\ell)=&125\cdot10^{8} \ell^8 + 2\cdot10^{8} \ell^6 r^2 - 100000 \ell^4 r^4 - 
 10400 \ell^2 r^6 - 3 r^8,\vspace{0.2cm}\\
 g_{2}(d,r,\ell)=&25 d^4 \ell^2 - 110000 d^2 \ell^4 - 4000000 \ell^6 - 1400 d^2 \ell^2 r^2 + 
 80000 \ell^4 r^2 -\\
& 3 d^2 r^4 + 1200 \ell^2 r^4.
\end{array}
$$
In addition the resultant of these two polynomials with respect to the variable $r$ has the factor
$$
g_{3}(d,\ell)=d^8 - 1240000 d^6 \ell^2 - 269440000 d^4 \ell^4 + 6144\cdot 10^{6} d^2 \ell^6 + 
 28672\cdot 10^{8} \ell^8.
$$
Consequently we will replace the problem of finding a zero of $N_{4}$, $N_{5}$ and $E_{2}$ by the problem of finding a zero of $g_{1}$, $g_{2}$ and $g_{3}$. Moreover we suppose that we have a common zero between each $g_{i}$ and $N_{3}$, $i=1,2,3$. With this supposition and computing the resultant between $g_{1}$ and $N_{3}$ we obtain a new polynomial $g_{4}$ in the variables $r$ and $d$ given by
$$
\begin{array}{rl}
g_{4}(r,d)=&1423828125 d^{16} + 25974000000 d^{14} r^2 - 47933388000000 d^{12} r^4 - \vspace{0.2cm}\\
 &478810245120000 d^{10} r^6 + 1711832243200000 d^8 r^8 - \vspace{0.2cm}\\
 &125591060860108800 d^6 r^{10} - 1334616713986048000 d^4 r^{12} + \vspace{0.2cm}\\
 &5191865574606503936 d^2 r^{14} - 4405603330689073152 r^{16}.
\end{array}
$$
However the resultant between $g_{3}$ and $g_{4}$ writes $Kr^{128}$, where $K$ is a positive constant. This polynomial cannot be zero since $r\geq0$. It means that when we consider each $g_{i}$ zero, $N_{3}$ must be positive or negative, exclusively. But taking the values $(r,\ell,d)=(r_{0},\ell_{0},1)$ we have $g_{i}\equiv0$, for each $i=1,2,3$ and $N_{3}$ is approximately $8738483/2500>0$, where $r_{0}$ and $\ell_{0}$ are roots of the polynomials
$$
\overline{r}(x)=-125 - 44400 x^2 + 184640 x^4 - 237568 x^6 + 86016 x^8,
$$
and
$$
\overline{l}(x)=1 - 1240000 x^2 - 269440000 x^4 + 6144000000 x^6 + 
 2867200000000 x^8.
$$
respectively. So $N_{3}$ is positive when $g_{i}\equiv0$, for each $i=1,2,3$ and then $D_{3}$ is negative. Hence using the factor $N_{5}$ of $D_{5}$ we cannot obtain 2 periodic solutions. The proof that cannot exists 2 periodic solutions for the case that we consider the factor $M_{5}$ of $D_{5}$ instead of $N_{5}$ is similar and we will omit it here.

\smallskip

We remark that we cannot prove analytically the non--existence of 5 periodic solutions for system \eqref{vdp_new} considering $m=0$, $n_{1}=n_{3}=n_{4}=1$ and $n_{2}>1$. Actually, using the algebraic manipulator Mathematica we obtained evidences that this number of periodic solutions cannot happen, but an analytical treatment is not trivial.

\smallskip

\subsection{Discussion of the results}\label{discussion}

Now we discuss some aspects of the averaging method presented in section \ref{averaging} in order to find periodic solutions in system \eqref{vdp}. In fact by using the averaging theory and the rescaling \eqref{rescaling} we shall see that we cannot prove that system \eqref{vdp} has periodic solutions except in the cases presented in Theorems from 1 to 8.

\smallskip

We start studying the possible values to the powers of $\varepsilon$ in system \eqref{vdp_new}. This is important because these powers play an important role in averaging theory described in section \ref{averaging} because they determine the terms of order 1 that we are interested in. Indeed we note that each one of the 5 different powers of $\varepsilon$ in system \eqref{vdp_new} must be non--negative. Therefore considering these powers as zero or positive we have 32 possible combinations of these 5 different powers of $\varepsilon$. Actually some of them are not algebraically possible. Table 1 exhibit only the possible case.

\smallskip

\begin{table}[h!]
\label{cases}
\normalsize
\centering
\renewcommand{\arraystretch}{1.2}
\begin{tabular}{|c|c|c|c|c|c|}
\hline
\multicolumn{6}{|c|}{Possible combination of powers of $\varepsilon$ for the system \eqref{F0_vdp}} \\
\hline
\hline
\multirow{16}{*}{$n_{1}=0$} & \multirow{8}{*}{$2m+n_{1}=0$} & \multirow{4}{*}{$2m+n_{2}=0$} &  \multirow{2}{*}{$4m+n_{3}=0$} & $-m+n_{4}=0$ & C1  \\
\cline{5-6}
 &  & &  & $-m+n_{4}>0$ & C2  \\
\cline{4-6}
 &  & & \multirow{2}{*}{$4m+n_{3}>0$} & $-m+n_{4}=0$ & C3  \\
\cline{5-6}
 &  & &  & $-m+n_{4}>0$ & C4  \\
\cline{3-6}
 & & \multirow{4}{*}{$2m+n_{2}>0$} & \multirow{2}{*}{$4m+n_{3}=0$} & $-m+n_{4}=0$ & C5  \\
\cline{5-6}
 &  & &  & $-m+n_{4}>0$ & C6  \\
\cline{4-6}
 &  & & \multirow{2}{*}{$4m+n_{3}>0$} & $-m+n_{4}=0$ & C7  \\
\cline{5-6}
 &  & &  & $-m+n_{4}>0$ & C8  \\
\cline{2-6}
 & \multirow{8}{*}{$2m+n_{1}>0$} & \multirow{4}{*}{$2m+n_{2}=0$} & \multicolumn{3}{|c|}{\multirow{4}{*}{$\emptyset$}}    \\  
 &  &  &  \multicolumn{3}{|c|}{}   \\ 
 &  &  &  \multicolumn{3}{|c|}{}   \\
 &  & &  \multicolumn{3}{|c|}{}   \\
\cline{3-6}
 & & \multirow{4}{*}{$2m+n_{2}>0$} & \multirow{2}{*}{$4m+n_{3}=0$} &  \multicolumn{2}{|c|}{\multirow{2}{*}{$\emptyset$}}  \\
 &  & &  &  \multicolumn{2}{|c|}{} \\
\cline{4-6}
 &  & & \multirow{2}{*}{$4m+n_{3}>0$} & $-m+n_{4}=0$ & C9  \\
\cline{5-6}
 &  & &  & $-m+n_{4}>0$ & C10  \\
\cline{1-6}
\multirow{16}{*}{$n_{1}>0$}  & \multirow{8}{*}{$2m+n_{1}=0$} &  \multicolumn{4}{|c|}{\multirow{8}{*}{$\emptyset$}} \\ 
 &   &  \multicolumn{4}{|c|}{}  \\
 &  &  \multicolumn{4}{|c|}{}  \\
 &  &  \multicolumn{4}{|c|}{}  \\
 &   &  \multicolumn{4}{|c|}{}  \\
 &  &  \multicolumn{4}{|c|}{}  \\
 &  &  \multicolumn{4}{|c|}{}  \\
 &  &  \multicolumn{4}{|c|}{}  \\
\cline{2-6}
 & \multirow{8}{*}{$2m+n_{1}>0$} & \multirow{4}{*}{$2m+n_{2}=0$} & \multirow{2}{*}{$4m+n_{3}=0$} & $-m+n_{4}=0$ & C11  \\
\cline{5-6}
 &  & &  & $-m+n_{4}>0$ & C12  \\
\cline{4-6}
 &  & & \multirow{2}{*}{$4m+n_{3}>0$} & $-m+n_{4}=0$ & C13  \\
\cline{5-6}
 &  & &  & $-m+n_{4}>0$ & C14  \\
\cline{3-6}
 & & \multirow{4}{*}{$2m+n_{2}>0$} & \multirow{2}{*}{$4m+n_{3}=0$} & $-m+n_{4}=0$ & C15  \\
\cline{5-6}
 &  & &  & $-m+n_{4}>0$ & C16  \\
\cline{4-6}
 &  & & \multirow{2}{*}{$4m+n_{3}>0$} & $-m+n_{4}=0$ & C17  \\
\cline{5-6}
 &  & &  & $-m+n_{4}>0$ & C18  \\
\cline{1-6}
\hline
\end{tabular}
\caption{Here, when appear $\emptyset$, this means that the corresponding cases cannot occur. The notation Ci, $i=1,\hdots,18$ enumerate each possible case.}
\end{table}

\smallskip

Observing Table 1 we see that in fact we have only 18 cases. Moreover in order to apply averaging theory, for which one of the 18 cases we must integrate the equations of the non--perturberd part of vector field of system \eqref{vdp_new}. This is not a simple task because the major part of these equations are nonlinear and non--autonomous. Indeed in each case from 1 to 8 and from 11 to 16 we could not integrate the equations even using the algebraic manipulators Mathematica or Maple. The non--integrable cases from Table 1 are listed above.
\begin{enumerate}
\item[{\it Caso 1:}] $F_{0}(t,x)=(y,-x+ry-ax^{3}-rx^{2}y-\ell x^{5}+d\cos\,t)$,
\item[{\it Caso 2:}] $F_{0}(t,x)=(y,-x+ry-ax^{3}-rx^{2}y-\ell x^{5})$,
\item[{\it Caso 3:}] $F_{0}(t,x)=(y,-x+ry-ax^{3}-rx^{2}y+d\cos\,t)$,
\item[{\it Caso 4:}] $F_{0}(t,x)=(y,-x+ry-ax^{3}-rx^{2}y)$,
\item[{\it Caso 5:}] $F_{0}(t,x)=(y,-x+ry-rx^{2}y-\ell x^{5}+d\cos\,t)$,
\item[{\it Caso 6:}] $F_{0}(t,x)=(y,-x+ry-rx^{2}y-\ell x^{5})$,
\item[{\it Caso 7:}] $F_{0}(t,x)=(y,-x+ry-rx^{2}y+d\cos\,t)$,
\item[{\it Caso 8:}] $F_{0}(t,x)=(y,-x+ry-rx^{2}y)$,
\item[{\it Caso 11:}] $F_{0}(t,x)=(y,-x-ax^{3}-\ell x^{5}+d\cos\,t)$,
\item[{\it Caso 12:}] $F_{0}(t,x)=(y,-x-ax^{3}-\ell x^{5})$,
\item[{\it Caso 13:}] $F_{0}(t,x)=(y,-x-ax^{3}+d\cos\,t)$,
\item[{\it Caso 14:}] $F_{0}(t,x)=(y,-x-ax^{3})$,
\item[{\it Caso 15:}] $F_{0}(t,x)=(y,-x-\ell x^{5}+d\cos\,t)$,
\item[{\it Caso 16:}] $F_{0}(t,x)=(y,-x-\ell x^{5})$.
\end{enumerate}

\smallskip

In particular we observe that cases 12, 14 and 16 turn the non--perturbed part of system \eqref{vdp_new} into a Hamiltonians system. In \cite{BL} the authors provide a method to apply averaging theory in planar systems when the system is Hamiltonian but the expressions become also to much complicate and again we could not integrate the equations in cases 12, 14 and 16.

\smallskip

In cases 9, 10 and 17 we can integrate the expressions but in these cases the hypotheses of averaging method do not apply. In these cases we have the following expressions for $F_{0}(t,{\bf x})$.

\begin{enumerate}
\item[{\it Caso 9:}] $F_{0}(t,{\bf x})=(y,-x+ry+d\cos\,t)$,
\item[{\it Caso 10:}] $F_{0}(t,{\bf x})=(y,-x+ry)$,
\item[{\it Caso 17:}] $F_{0}(t,{\bf x})=(y,-x+d\cos\,t)$.
\end{enumerate}

Solution $x(t)$ in case 9 is
$$
C^{-}e^{\dfrac{1}{2}t\left(r-\sqrt{-4+r^{2}}\right)}+C^{+}e^{\dfrac{1}{2}t\left(r+\sqrt{-4+r^{2}}\right)}-(d/r)\sin\,t,
$$
where
$$
C^{\mp}=\dfrac{\mp 2d+r\left(\left(  \pm r+\sqrt{-4+r^{2}} \right)x_{0}\mp2y_{0}\right)  }{2r\sqrt{-4+r^{2}}}.
$$
In order that $x(t)$ be periodic we must have $C^{-}$ and $C^{+}$ equal to zero. But these conditions are verified only if we take $x_{0}=0$ and $y_{0}=-d/r$. So $x_{0}$ and $y_{0}$ are fixed and then the unperturbed system $\dot{x}=F_{0}(t,x)$ has no submanifold of periodic solutions. Consequently we can not apply averaging theory in this case.

\smallskip

On the other hand, in case 10 solution $x(t)$ writes
$$
C^{-}e^{\dfrac{1}{2}t\left(r-\sqrt{-4+r^{2}}\right)}+C^{+}e^{\dfrac{1}{2}t\left(r+\sqrt{-4+r^{2}}\right)}-(d/r)\sin\,t,
$$
where
$$
C^{\mp}=\dfrac{\left(  \pm r+\sqrt{-4+r^{2}} \right)x_{0}\mp2y_{0}  }{2r\sqrt{-4+r^{2}}}.
$$
Again we need to chose $C^{-}$ and $C^{+}$ equal to zero but this happen only if $(x_{0},y_{0})=(0,0)$. So we cannot apply averaging theory because $(0,0)$ is the equilibrium point of the non--perturbed part of system \eqref{vdp_new} in case 10.

\smallskip

Finally in case 17 solutions $x(t)$ and $y(t)$ have the form
$$
\begin{array}{l}
x(t)=x_{0}\cos\,t+y_{0}\sin\,t+\dfrac{d}{2}t\cos\,t,\vspace{0.2cm}\\
y(t)=y_{0}\cos\,t-x_{0}\sin\,t+\dfrac{d}{2}(t\cos\,t+\sin\,t).
\end{array}
$$
It is immediate to note that these solutions are non--periodic because $d$ is positive.

\smallskip

\smallskip

Case 18 provide positive results and Theorems from 1 to 8 are based in this case. In fact we observe that in case 18 the vector field of system \eqref{vdp_new} is
$$
(y,-x+\varepsilon^{n_{1}}ry-\varepsilon^{2m+n_{1}}rx^{2}y-\varepsilon^{2m+n_{2}}ax^{3}-\varepsilon^{4m+n_{3}}\ell x^{5}+\varepsilon^{-m+n_{4}}d\cos\,t)^{T}.
$$

\smallskip

As we commented before the periodic solutions of our problem correspond to the simple zeros of function \eqref{average_function}. So in order to calculate the zeros of the function $f({\bf z})$ given in \eqref{average_function} we must determine $F_{1}(t,{\bf x})$, where $F_{1}(t,{\bf x})$ is determinated by the terms of order 1 on $\varepsilon$ of the vector field of system \eqref{vdp_new}. However we see that the expression of $F_{1}(t,{\bf x})$ depends on the values $m$ and $n_{i}$, $i=1,2,3,4$. In fact, if $m$ is positive and observing the conditions \eqref{conditions} the only terms of the vector field of order 1 on $\varepsilon$ can be generate by powers $n_{1}$ and $-m+n_{4}$, where these two values can be 1 or greater than 1. This implies four possibilities. On the other hand, if $m$ is zero the powers of $\varepsilon$ in the vector field of system \eqref{vdp_new} depend on $n_{1}$, $n_{2}$, $n_{3}$ and $n_{4}$. Again, each one of these powers can be 1 one greater than 1. So, if $m$ is zero we have 16 possibilities for $F_{1}$. Then case 18 has 20 subcases corresponding to the different possibilities of $F_{1}(t,{\bf x})$. These subcases are presented in Table 2.

\begin{table}[htb]
\label{subcases}
\normalsize
\centering
\renewcommand{\arraystretch}{1.2}
\begin{tabular}{|c|c|c|}
\hline
Subc. & Conditions & Second coordinate of $F_{1}(t,{\bf x})$ \\
\hline
\hline
1 & $m=0$, $n_{i}=1$, $i=1,2,3,4$ & $ry-ax^3-rx^{2}y-\ell x^{5}+d\cos\,t$  \\
\hline
2 & $m=0$, $n_{i}=1$, $i=1,2,3$, $n_{4}>1$ & $ry-ax^3-rx^{2}y-\ell x^{5}$  \\
\hline
3 & $m=0$, $n_{i}=1$, $i=1,2,4$, $n_{3}>1$ & $ry-ax^3-rx^{2}y+d\cos\,t$  \\
\hline
4 & $m=0$, $n_{i}=1$, $i=1,3,4$, $n_{2}>1$ & $ry-rx^{2}y-\ell x^{5}+d\cos\,t$  \\
\hline
5 & $m=0$, $n_{i}=1$, $i=2,3,4$, $n_{1}>1$ & $-ax^3-\ell x^{5}+d\cos\,t$  \\
\hline
6 & $m=0$, $n_{1}=n_{2}=1$, $n_{3},n_{4}>1$ & $ry-ax^3-rx^{2}y$  \\
\hline
7 & $m=0$, $n_{1}=n_{3}=1$, $n_{2},n_{4}>1$ & $ry-rx^{2}y-\ell x^{5}$  \\
\hline
8 & $m=0$, $n_{1}=n_{4}=1$, $n_{2},n_{3}>1$ & $ry-rx^{2}y+d\cos\,t$  \\
\hline
9 & $m=0$, $n_{2}=n_{3}=1$, $n_{1},n_{4}>1$ & $-ax^3-\ell x^{5}$  \\
\hline
10 & $m=0$, $n_{2}=n_{4}=1$, $n_{1},n_{3}>1$ & $-ax^3+d\cos\,t$  \\
\hline
11 & $m=0$, $n_{3}=n_{4}=1$, $n_{1},n_{2}>1$ & $-\ell x^{5}+d\cos\,t$  \\
\hline
12 & $m=0$, $n_{1}=1$, $n_{i}>1$, $i=2,3,4$ & $ry-rx^{2}y$  \\
\hline
13 & $m=0$, $n_{2}=1$, $n_{i}>1$, $i=1,3,4$ & $-ax^3$  \\
\hline
14 & $m=0$, $n_{3}=1$, $n_{i}>1$, $i=1,2,4$ & $-\ell x^{5}$  \\
\hline
15 & $m>0$, $n_{4}=1$, $n_{i}>1$, $i=1,2,3$ & $d\cos\,t$  \\
\hline
16 & $m>0$, $n_{i}>1$, $i=1,2,3,4$ & $0$  \\
\hline
17 & $m>0$, $n_{1}=1$, $-m+n_{4}=1$ & $ry+d\cos\,t$  \\
\hline
18 & $m>0$, $n_{1}=1$, $-m+n_{4}>1$ & $ry$  \\
\hline
19 & $m>0$, $n_{1}>1$, $-m+n_{4}=1$ & $d\cos\,t$  \\
\hline
20 & $m>0$, $n_{1}>1$, $-m+n_{4}>1$ & $0$  \\
\hline
\end{tabular}
\caption{Possible expressions for $F_{1}(t,{\bf x})$ when $m^2+n_{2}^2$, $m^2+n_{3}^2$, $n_{4}-m$ and $n_{1}$ are positives. We exhibit only the second coordinate of $F_{1}(t,{\bf x})$ because the first one has no terms depending on $\varepsilon$.}
\end{table}

\smallskip

Theorems from 1 to 8 correspond to subcases 8, 10, 11, 17, 5, 3, 1 and 4, respectively. These are the only subcases of case 18 that the averaging method provide positive results. Indeed, in cases 15, 16, 19 and 20 function \eqref{average_function} do not have zeros. On the other hand, in cases 2, 6, 7, 13, 14, 18 the only zero of function \eqref{average_function} is $(x_{0},y_{0})=(0,0)$, that correspond to the equilibrium point of the system, and consequently in these subcases the system does not have periodic solutions. In cases 9 and 12 function \eqref{average_function} has real zeros different from $(0,0)$ but they are non--isolated and then the Jacobian of function \eqref{average_function} at the zero is zero, and consequently the averaging theory cannot be applied.

\smallskip

\section{The Averaging Theory for Periodic Solutions}\label{averaging}

\smallskip

Now we present the basic results on the averaging theory of first order that we need to proving our results.

\smallskip

Consider the problem of bifurcation of $T$--periodic solutions from differential systems of the form
\begin{equation}\label{perturbed}
\dot{{\bf x}}=F_{0}(t,{\bf x})+\varepsilon F_{1}(t,{\bf x})+\varepsilon^{2} R(t,{\bf x},\varepsilon),
\end{equation}
with $\varepsilon=0$ to $\varepsilon\neq0$ sufficiently small. Here the functions $F_{0},F_{1}:\mathbb{R}\times\Omega\rightarrow \mathbb{R}^{n}$ and $R:\mathbb{R}\times\Omega\times(-\varepsilon_{f},\varepsilon_{f})\rightarrow \mathbb{R}^{n}$ are $\mathcal{C}^{2}$, $T$--periodics in the first variable and $\Omega$ is an open subset of $\mathbb{R}^{n}$. One of the main assumption is that the unperturbed system
\begin{equation}\label{unperturberd}
\dot{{\bf x}}=F_{0}(t,{\bf x}),
\end{equation}
has a manifold of periodic solutions. A solution of this problem is given using averaging theory.

\smallskip

Indeed, assume that there is an open set $V$ with $\overline{V}\subset D\subset\Omega$ and such that for each ${\bf z}\in\overline{V}$, ${\bf x}(\cdot,{\bf z},0)$ is $T-$periodic, where ${\bf x}(\cdot,{\bf z},0)$ is the solution of the unperturbed system \eqref{unperturberd} with ${\bf x}(0)={\bf z}$. As answer to the problem of bifurcation of $T-$periodic solutions from ${\bf x}(\cdot,{\bf z},0)$ is given in the following Theorem.



\begin{theorem}\label{averaging1order}
We assume that there exists an open set $V$ with $\overline{V}\subset D$ and such that for each ${\bf z}\in\overline{V}$, ${\bf x}(\cdot,{\bf z},0)$ is $T-$periodic and consider the function $f:\overline{V}\rightarrow \mathbb{R}^{n}$ given by
$$
f({\bf z})=\displaystyle\int_{0}^{T}Y^{-1}(t,{\bf z})F_{1}(t,{\bf x}(t,{\bf z},0))dt.
$$
Then the following statements hold.
\begin{enumerate}
\item[(a)] If there exists $a\in V$ with $f(a)=0$ and $\det((\partial f/\partial {\bf z})(a))\neq0$, then there exists a $T-$periodic solution $\varphi(\cdot,\varepsilon)$ of system \eqref{perturbed} such that $\varphi(0,\varepsilon)\to a$ as $\varepsilon\to 0$.
\item[(b)] The type of stability of the periodic solution $\varphi(\cdot,\varepsilon)$ is given by the eigenvalues of the Jacobian matrix $M=((\partial f/\partial {\bf z})(a))$.
\end{enumerate}

\end{theorem}
For a proof of Theorem \ref{averaging1order}(a) see Corollary 1 of \cite{BFL}.

In fact, the result of Theorem \ref{averaging1order} is a classical result due to Malkin \cite{MALKIN} and Roseau \cite{ROSEAU}. For a shorter proof of Theorem \ref{averaging1order}$(a)$, see \cite{BFL}.

\smallskip

For Additional information on averaging theory see the book \cite{SVM}.

\section*{Acknowledgments}

The first author is supported by the FAPESP-BRAZIL   grants 2010/18015-6 and 2012/05635-1. The second author is partially supported by the grants MINECO/FEDER MTM 2008--03437, AGAUR 2009SGR 410, ICREA  \linebreak Academia and FP7  PEOPLE-2012-IRSES-316338 and 318999 and CAPES--MECD grant PHB-2009-0025-PC.

\appendix
\section{Root Classification for a Quintic Polynomial}\label{appendix}

In this section we present a brief summary of the results about the number and multiplicities of the real/complex roots for a quintic polynomial with arbitrary coefficients presented in \cite{ROOTS}. Indeed, consider the polynomial $P(x)=x^{5}+px^{3}+qx^{2}+ux+v$. So, the following table gives the number of real and complex roots and multiplicities roots of $P(x)$ in all cases.

$$
\begin{array}{cll}
(1) &D_{5}>0\;\wedge\;D_{4}>0\;\wedge\;D_{3}>0\;\wedge\;D_{2}>0 & \{1,1,1,1,1\}\\
(2) &D_{5}>0\;\wedge\;(D_{4}\leq0\;\vee\;D_{3}\leq0\;\vee\;D_{2}\leq0)&\{1\}\\
(3) &D_{5}<0&\{1,1,1\}\\
(4) &D_{5}=0\;\wedge\;D_{4}>0&\{2,1,1,1\}\\
(5) &D_{5}=0\;\wedge\;D_{4}<0&\{2,1\}\\
(6) &D_{5}=0\;\wedge\;D_{4}=0\;\wedge\;D_{3}>0\;\wedge\;E_{2}\neq 0&\{2,2,1\}\\
(7) &D_{5}=0\;\wedge\;D_{4}=0\;\wedge\;D_{3}>0\;\wedge\;E_{2}=0&\{3,1,1\}\\
(8) &D_{5}=0\;\wedge\;D_{4}=0\;\wedge\;D_{3}<0\;\wedge\;E_{2}\neq0&\{1\}\\
(9) &D_{5}=0\;\wedge\;D_{4}=0\;\wedge\;D_{3}<0\;\wedge\;E_{2}=0&\{3\}\\
(10) &D_{5}=0\;\wedge\;D_{4}=0\;\wedge\;D_{3}=0\;\wedge\;D_{2}\neq0\;\wedge\;F_{2}\neq0&\{3,2\}\\
(11) &D_{5}=0\;\wedge\;D_{4}=0\;\wedge\;D_{3}=0\;\wedge\;D_{2}\neq0\;\wedge\;F_{2}=0&\{4,1\}\\
(12) &D_{5}=0\;\wedge\;D_{4}=0\;\wedge\;D_{3}=0\;\wedge\;D_{2}=0&\{5\}
\end{array}
$$

where
$$
\begin{array}{rl}
D_{2}=&-p,\vspace{0.2cm}\\
D_{3}=&-12 p^3 - 45 q^2 + 40 p u,\vspace{0.2cm}\\
D_{4}=&-4 p^3 q^2 - 27 q^4 + 12 p^4 u + 117 p q^2 u - 88 p^2 u^2 + 160 u^3 - 
 40 p^2 q v - \vspace{0.2cm}\\
 &300 q u v + 125 p v^2,\vspace{0.2cm}\\
D_{5}=&-4 p^3 q^2 u^2 - 27 q^4 u^2 + 16 p^4 u^3 + 144 p q^2 u^3 - 
 128 p^2 u^4 + 256 u^5 + \vspace{0.2cm}\\
 &16 p^3 q^3 v + 108 q^5 v - 72 p^4 q u v - 
 630 p q^3 u v + 560 p^2 q u^2 v - 1600 q u^3 v +\vspace{0.2cm}\\
 & 108 p^5 v^2 + 
 825 p^2 q^2 v^2 - 900 p^3 u v^2 + 2250 q^2 u v^2 + 2000 p u^2 v^2 -\vspace{0.2cm}\\
 & 3750 p q v^3 + 3125 v^4,\vspace{0.2cm}\\
E_{2}=&16 p^4 q^2 - 48 p^5 u + 60 p^2 q^2 u + 160 p^3 u^2 + 900 q^2 u^2 - 
 1100 p^3 q v -\vspace{0.2cm}\\
 & 3375 q^3 v + 1500 p q u v + 625 p^2 v^2,\vspace{0.2cm}\\
F_{2}=&3q^{2}-8pu.
\end{array}
$$

The polynomials $D_{i}$, $i=2,3,4,5$, $E_{2}$ and $F_{2}$ form a discriminant system which is sufficient for the classification of roots of the polynomial $P(x)$, which is described by the right column of the table. For instance, $\{1,1,1\}$ means three real simple roots and a pair of complex roots, and $\{3,1,1\}$ means a real triple root plus two real simple roots.


\begin{thebibliography}{99}

\bibitem{AL} {\sc W. W. Adams and P. Loustaunau},
{\it An Introduction to Gr\"{o}bner Bases}, American Mathematical
Society, Graduate Studies in Mathematics, Vol. {\bf 3}, 1994.

\bibitem{GUCKII} {\sc K. Bold, C. Edwards, J. Guckenheimer, S. Guharay, K. Hoffman, J. Hubbard, R. Oliva and W. Weckesser},
{\it The forced Van der Pol equation II: Canards in the reduced system}, SIAM J. Appl. Dyn. Syst. {\bf 2} (2003), 570--608.

\bibitem{BFL} {\sc A. Buic\v{a}, J. P. Fran\c{c}oise and J. Llibre},
{\it Periodic solutions of nonlinear periodic differential systems with a small parameter}, Comm. on Pure and Appl. Anal. {\bf 6} (2007), 103--111.

\bibitem{BL} {\sc A. Buic\v{a} and J. Llibre},
{\it Averaging methods for finding periodic orbits via Brouwer degree}, Bull. Sci. Math. {\bf 128} (2004), 7--22.

\bibitem{CJ} {\sc A. Chen and G. Jiang},
{\it Periodic solution of the Duffing-Van der Pol oscillator by homotopy perturbation method}, Int. J. Comput. Math {\bf  87} (2010), 2688--2696.

\bibitem{CL} {\sc Y. Chen and J. Liu},
{\it Uniformly valid solution of limit cycle of the Duffing--van der Pol equation}, Mech. Res. Comm. {\bf 36} (2009), 845--850.

\bibitem{CJN} {\sc W. O. Criminale, T. L.  Jackson and P. W. Nelson},
{\it Limit cycle-strange attractor competition}, Stud. Appl. Math. {\bf 112}, (2004), 133--160.

\bibitem{EH} {\sc C. Egami and N. Hirano},
{\it Periodic solutions for forced van der Pol type equations}, Mathematical economics {\bf 1264}, (2001), 159--172.

\bibitem{EL} {\sc R. D. Euzébio and J. Llibre},
{\it Periodic Solutions of \textit{El Ni\~{n}o} Model through the Vallis Differential System}, preprint.

\bibitem{DIDATIC} {\sc T. H. Fay},
{\it The forced Van der Pol equation}, Int. Journal of Math. Educ. in Sci. and Tech. {\bf 40}(5) (2009), 669--677.

\bibitem{GUCKI} {\sc J. Guckenheimer, K. Hoffman and W. Weckesser},
{\it The forced Van der Pol equation I: The slow flow and its bifurcations}, SIAM J. Appl. Dyn. Syst. {\bf 2} (2003), 1--35.

\bibitem{DV7} {\sc F. M. M. Kakmeni, S. Bowong, C. Tchawoua and E. Kaptouom},
{\it Strange attractors and chaos control in a Duffing--Van der Pol oscillator with two external periodic forces}, J. Sound Vib. {\bf 227}(4--5) (2004), 783--799.

\bibitem{PRINCIPAL20} {\sc F. M. M. Kakmeni, S. Bowong, C. Tchawoua and E. Kaptouom},
{\it Chaos control and synchronization of a $\phi^6$-Van der Pol oscillator}, Phys. Lett. A {\bf 322}(5--6) (2004), 305--323.

\bibitem{LEUNG} {\sc H. K. Leung},
{\it Synchronization dynamics of coupled van der Pol systems}, Physica A. {\bf 321} (2003), 248--255.

\bibitem{L2} {\sc H. Li},
{\it Gr\"{o}bner Bases in Ring Theory},  World Scientific Publishing,
2011.

\bibitem{LXZ} {\sc Z. Li, W. Xu and X. Zhang},
{\it Analysis of chaotic behavior in the extended Duffing--Van der Pol system subject to additive non-symmetry biharmonical excitation}, Appl. Math. Comput. {\bf 183} (2006), 858--871.

\bibitem{PRINCIPAL} {\sc D. Liu and H. Yamaura},
{\it Chaos control of a $\phi^6$-Van der Pol oscillator driven by external excitation}, Nonlinear Dynamics {\bf 68}(1--2) (2012), 95--105.

\bibitem{LV}
{\sc J. Llibre and C. Vidal}, {\it Periodic solutions of a periodic FitzHugh-Nagumo differential system}, preprint.

\bibitem{GENERAL} {\sc M. Ma, J. Zhou and J. Cai},
{\it Practical synchronization of second-order non--autonomous systems with parameter mismatch and its applications}, Nonlinear Dyn. {\bf 69} (2012), 1285--1292.

\bibitem{MALKIN} {\sc I. G. Malkin},
{\it Some problems of the theory of nonlinear oscillations}, (Russian) Gosudarstv. Izdat. Tehn.-Teor. Lit., Moscow, 1956.

\bibitem{ROSEAU} {\sc M. Roseau},
{\it Vibrations non lin\'{e}aires et th\'{e}orie de la stabilit\'{e}}, (French) Springer Tracts in Natural Philosophy, Vol. 8 Springer--Verlag, Berlin--New York, 1966.

\bibitem{PRINCIPAL21} {\sc H. Salarieh and A. Alasty},
{\it Control of stochastic chaos using sliding mode method}, J. Comput. Appl. Math. {\bf 225} (2009), 135--145.

\bibitem{SVM} {\sc J. A. Sanders, F. Verhulst and J. Murdock},
{\it Averaging method in nonlinear dynamical systems}, Appl. Math.
Sci., vol. {\bf 59}, Springer, New York, 2007.

\bibitem{COMPLEX9} {\sc M. S. Siewe, F. M. Kakmeni and C. Tchawoua},
{\it Resonant os- cillation and homoclinic bifurcation in a$\phi^6$-Van der Pol oscillator}, Chaos, Solitons and Fractals {\bf 21}(4) (2004), 841--853.

\bibitem{DV5} {\sc W. Szemplinska-Stupnicka and J. Rudowski},
{\it Neimark bifurcation almost-periodicity and chaos in the forced Van der Pol-duffing system in the neighbourhood of the principal resonance}, Phys. Lett. A {\bf 192}(2--4) (1994), 201--206.

\bibitem{COMPLEX8} {\sc R. Tchoukuegno, B. R. N. Nbendjo and P. Woafo},
{\it Resonant oscillations and fractal basin boundaries of a particle in a $\phi^6$-potential}, Physica A {\bf 304}(3--4) (2002), 362--378.

\bibitem{DV4} {\sc Y. Ueda and N. Akamatsu},
{\it Chaotically transitional phenomena in the forced negative-resistance oscillator}, IEEE Trans. Circuits Syst. {\bf 28}(3) (1981), 217--224.

\bibitem{DV6} {\sc A. Venkatesan and M. Lakshmanan},
{\it Bifurcation and chaos in the double-well Duffing--Van der Pol oscillator: numerical and analytical studies.}, Phys. Rev. E {\bf 56}(6) (1997), 6321--6330.

\bibitem{ROOTS} {\sc L. Yang},
{\it Recent Advances on Determining the Number of Real Roots of Parametric Polynomials}, J. Symbolic Computation {\bf 28} (1999), 225--242.

\bibitem{YXS} {\sc X. Yang, W. Xu and Z. Sun},
{\it Effect of Gaussian white noise on the dynamical behaviours of an extended Duffing-Van der Pol oscillator}, Internat. J. Bifur. Chaos Appl. Sci. Engrg. 16 {\bf 9} (2006), 2587--2600.

\bibitem{YXS2} {\sc X. Yang, W. Xu and Z. Sun},
{\it Effect of bounded noise on the chaotic motion of a Duffing Van der pol oscillator in a $\phi^{6}$ oscillator}, Chaos, Solitons and Fractals 27 {\bf 3} (2006), 778--788.


















\end{thebibliography}
\end{document}